\documentclass{amsart}
    \usepackage[margin=0.7in]{geometry}
    \usepackage[parfill]{parskip}
    \usepackage[utf8]{inputenc}
    \usepackage{hyperref}
    \usepackage{tikz}
    \usepackage{xcolor}
    \usepackage{pgfplots}
    \usepackage{subcaption}
    \usepackage{filecontents}
    \usetikzlibrary{patterns}
    \usepackage{slashbox}

    \usepackage{amsmath,amssymb,amsfonts,amsthm}

\DeclareMathOperator*{\argmax}{arg\,max}
    
\author{Denys Dragunov}
\address[Denys Dragunov]{Kyiv, Ukraine}
\email{dragunovdenis@gmail.com}
\dedicatory{To my dear daughter Viktoriia.}

\title{A transformation-based approach for solving stiff two-point boundary value problems}
\begin{document}

\begin{abstract}
A new approach for solving stiff boundary value problems for systems of ordinary differential equations is presented. Its idea essentially generalizes and extends that from \cite{Makarov_Dragunov_2019}. The approach can be viewed as a methodology framework that allows to enhance "stiffness resistance" capabilities of pretty much all the known numerical methods for solving two-point BVPs. The latter is demonstrated on the example of the {\it trapezoidal scheme} with the corresponding C++ source code available at \url{https://github.com/imathsoft/MathSoftDevelopment}. Results of numerical experiments are provided to support the theoretical conclusions. 
\end{abstract}
\maketitle
\large

\section{Introduction}\label{introduction_section}
Stiffness, as a property of a boundary value problem (BVP), can manifest itself in a number of different ways. It might be that the variety of forms the stiffness can take is eventually responsible for the fact that until the recent times there were no strict formal definition of the phenomenon \cite{Fifty_years_of_stiffness}. In different sources one can find rather informal and "pragmatic" definitions of stiff problems, like those "for which explicit methods don't work" \cite{Hairer_Wanner_Stiff}.
Aiming to be practically oriented and concrete, the current paper employs its own (informal and by no means complete or general) definition of a stiff two-point BVP through the particular behaviours of the problem's solution. Namely, throughout this paper, saying that a two-point BVP is stiff we mean that its solution has either {\it boundary layers} or {\it shock layers}. By a {\it boundary layer} we mean a boundary-adjacent narrow interval within the solution's domain, where the absolute value of the solution's derivative rapidly increases to "extremely" high magnitudes. A similar interval which is not incident to either of the two boundary points will be referred to as a {\it shock layer}.

Possessing broad practical applications, stiff BVPs have always been under scrutiny of the numerical analysis community and, as a result, there is a number of quite impressive methods and software packages available to deal with this type of problems (see, for example, \cite{Hairer_Wanner_Stiff}, \cite{Cash_bvpSolve_in_R_2014}, \cite{Cash_bvpSolve_test_problems_2010}). Almost all the numerical methods for solving stiff BVPs are, to a greater or lesser extent, focused on the "construction of a mesh on which all features of the solution are locally smooth" \cite{Lee_June_Yub_and_Greengard_Leslie_1997_SIAM}. The latter general strategy finds its implementation in a variety of ways, for example, by means of so called {\it monitor functions} \cite{Lee_June_Yub_and_Greengard_Leslie_1997_SIAM}, \cite{Wright_Cash_Moore_1994} as well as through the properly chosen smooth transformations applied to either the independent variable \cite{Kreiss_Nichols_Brown} or to the unknown function (the solution) \cite{Chang20103043}, \cite{Chang20103303}, \cite{Gen_sol_of_TP}. 

Following the forementioned "mesh-adjustment" paradigm, in \cite{Makarov_Dragunov_2019} yet another approach for dealing with stiff BVPs was suggested. Its idea naturally follows from an observation that the computational complexity (i.e. the stiffness) induced by the {\it boundary} and {shock layers} is largely due to the fundamental difficulties that most approximation methods experience when being applied to functions with large {\it moduli of smoothness } \cite{Moduli_of_smoothness_1988}. On the other hand, within the boundary/shock layers the solution is, obviously, monotone, and thus, its inverse is well defined and variates pretty moderately (being "almost" constant, since the derivative of an inverse function is equal to the reciprocal of the derivative of the original function)! With this in mind, paper \cite{Makarov_Dragunov_2019}, which is exclusively focused on the BVPs for the second order ordinary differential equations of the form
\begin{equation}\label{Second_order_BVP}
	u^{\prime\prime}(t) = N(u^{\prime}(t), u(t), t),\; t\in [a,b,]
\end{equation}
$$
u(a) = u_{a},\; u(b) = u_{b},
$$
proposes to switch to the equation for the "inverse solution" 
\begin{equation}\label{Inverse_second_order_BVP}
	t^{\prime\prime}(u) = -N(1/t^{\prime}(u), u, t(u))\left(t^{\prime}(u)\right)^{3}, \; t(\cdot) \stackrel{def}{=} u^{-1}(\cdot)
\end{equation}
within the boundary/shock layers. The approach demonstrated remarkable results on different test problems from \cite{Cash_bvpSolve_test_problems_2010}, including the well known Troesch's problem \cite{Troesch1976279}. Later the results from \cite{Makarov_Dragunov_2019} have been extended in \cite{2018arXiv181209498M}, mainly, in the part of theoretical justification and error analysis. 

It is worth mentioning, however, that the approach, in the form it was presented in \cite{Makarov_Dragunov_2019}, looks rather artificial and restricted. First of all, this is due to it being introduced specifically for the case of the second order ODEs and, second of all, --- due to rather sophisticated and unconventional implementation framework suggested by the authors. The present paper is aimed to address both this issues. Namely, the {\it transformation-based} approach, presented below, is a natural generalization of the one from \cite{Makarov_Dragunov_2019} for the case of systems of ODEs and, additionally, the particular implementation of the new technique, discussed here, is based on one of the most "canonical" methods for solving BVPs: the {\it trapezoidal scheme} \cite{Ascher_1988}. The latter is an attempt to view the suggested approach as a "stiffness-resistant enhancement" that can be applied to pretty much all the known methods for solving two-point BVPs.            

The paper is organized as follows. Section \ref{transformations_section} introduces the, so called, "swap" and "flip" transformations and discusses their stiffness-suppressing properties. Section \ref{Implementation_aspect_section} is focused on the implementation practicalities of the suggested approach, in particular, for the case when the trapezoidal scheme is chosen as the "base method". The results of numerical experiments demonstrating the "stiffness-resistance" capabilities of the transformation-based approach are presented in section \ref{Numerical_experiments_section}. Finally, some thoughts summarizing the article can be found in section \ref{Conclusions_section}. 

\section{The transformations}\label{transformations_section}
Let $\mathbb{S}_{n}$ denote a set of all systems of $n$ first order ordinary differential equations. Let's consider an arbitrary element $F$ from $\mathbb{S}_{n},$ which obviously can be expressed in the following form 
\begin{equation}\label{system}
F: \;\mathbf{u}^{\prime}(t) = \mathbf{F}(\mathbf{u}(t), t),
\end{equation}
where 
$$\mathbf{u}(t) = 
\left[\begin{array}{cccc}
 u_{1}(t), \ldots, u_{n}(t)
\end{array}\right]^{T},$$
$$\mathbf{F}(\mathbf{u}(t), t) = 
\left[\begin{array}{ccccc}F_{1}\left(\mathbf{u}(t), t\right), \ldots, F_{n}\left(\mathbf{u}(t), t\right) \end{array}\right]^{T},$$
$$F_{j}\left(\mathbf{u}(t), t\right) = F_{j}(u_{1}(t), \ldots, u_{n}(t), t),\; j\in \overline{1, n}.$$

For any $k\in \overline{1,n},$ we define a {\it $k$-swap} operator $\mathcal{SP}_{k}$ as a mapping $\mathbb{S}_{n} \rightarrow \mathbb{S}_{n}$ acting in the following way
\begin{equation}\label{system_transformed_swap}
G \stackrel{def}{=} \mathcal{SP}_{k}\left(F\right): \; \mathbf{v}^{\prime}(u) = \mathbf{G}(\mathbf{v}(u), u),
\end{equation} 	
where
\begin{equation}\label{vector_function_v}
\mathbf{v}(t) = 
\left[\begin{array}{cccc}
 v_{1}(u), \ldots, v_{n}(u)
\end{array}\right]^{T},
\end{equation}
$$\mathbf{G}(\mathbf{v}(u), u) = 
\left[\begin{array}{ccccc}G_{1}\left(\mathbf{v}(u), u\right), \ldots, G_{n}\left(\mathbf{v}(u), u\right) \end{array}\right]^{T},$$
$$G_{j}(\mathbf{v}(u), u) = \frac{F_{j}(v_{1}(u),\ldots, v_{k-1}(u), u, v_{k+1}(u), \ldots v_{n}(u), v_{k}(u))}{F_{k}(v_{1}(u),\ldots, v_{k-1}(u), u, v_{k+1}(u), \ldots v_{n}(u), v_{k}(u))}, \; j\in \overline{1, n},\; j \neq k,$$
$$G_{k}(\mathbf{v}(u), u) = \frac{1}{F_{k}(v_{1}(u),\ldots, v_{k-1}(u), u, v_{k+1}(u), \ldots v_{n}(u), v_{k}(u))}.$$

Let vector-function $\mathbf{u}(t),\; t\in [a, b]$ be a solution to system \eqref{system} (satisfying some initial or boundary conditions) and let function $u_{k}(t)$ be strictly monotone on some subinterval  $$(t_{0}, t_{1}) \subset (a, b).$$ It is easy to verify that, in this case, vector-function $\mathbf{v}(u)$ \eqref{vector_function_v}, defined as
\begin{equation}\label{v_definition}
v_{i}(u) = u_{i}(v_{k}(u)),\; i\in\overline{1,n},\; i \neq k,\; v_{k}(u) = u^{-1}_{k}(u),
\end{equation}
satisfies system \eqref{system_transformed_swap} $\forall u \in (u_{k}(t_{0}), u_{k}(t_{1})).$
As one can notice from the example above, the {\it $k$-swap} operator results in the unknown function $u_{k}(\cdot)$ and the independent variable $t$ being "swapped" (hence, the term), in the sense that in system \eqref{system_transformed_swap} the former becomes the new independent variable $u$, whereas the latter turns into a new unknown function $v_{k}(u)$. A practical application of such a "trick" becomes clear from the reasoning below. 

Let the forementioned interval $(t_{0}, t_{1})$ be a boundary or shock layer for the solution $\mathbf{u}(t).$ Then, by definition,
$$\|\mathbf{F}(\mathbf{u}(t), t)\| \gg 1,\; t\in (t_{0}, t_{1})$$
and at least for a single $k \in \overline{1, n},$ the inequality
$$|F_{k}(u_{1}(t),\ldots, u_{n}(t), t)| \gg 1,\; t\in (t_{0}, t_{1})$$
holds true. Provided that the interval is narrow enough, without loss of generality, we may assume that   
$$|F_{k}(u_{1}(t),\ldots, u_{n}(t), t)| \geq |F_{l}(u_{1}(t),\ldots, u_{n}(t), t)|, \; \forall l \in \overline{1, n}.$$ The latter inequality obviously implies that function $\mathbf{v}(u)$ \eqref{vector_function_v}, \eqref{v_definition}, satisfying system \eqref{system_transformed_swap}, can't have neither boundary nor shock layers on the interval $$\left(\min\{u_{k}(t_{0}), u_{k}(t_{1})\}, \max\{u_{k}(t_{0}), u_{k}(t_{1})\}\right),$$ which corresponds to interval $(t_{0}, t_{1})$ in terms of system \eqref{system}. 

To illustrate how the introduced "swap" transformation can be applied in practice, let us consider the well known Troesch's problem \cite{Troesch1976279}:
\begin{equation}\label{Troesch_problem}
\left[\begin{array}{l}
u_{1}^{\prime}(t) \\
u_{2}^{\prime}(t)
\end{array}\right] = \mathbf{T}_{0}(u_{1}(t), u_{2}(t), \lambda) \stackrel{def}{=} \left[\begin{array}{l}
u_{2}(t) \\
\lambda \sinh\left(\lambda u_{1}(t)\right)
\end{array}\right],\; t\in [0, 1], \; u_{1}(0) = 0,\; u_{1}(1) = 1.
\end{equation}
The problem has a boundary layer near the right boundary point $t=1$ so that $$u^{\prime}_{i}(t) \gg 1, \; \forall t \in (\varepsilon, 1),\; i=1,2, $$
provided that $\varepsilon < 1$ is close enough to $1$ and parameter $\lambda$ is big enough. The stiffness of the problem increases along with parameter $\lambda,$ which, in effect, amounts to the shrinkage of interval $(\varepsilon, 1)$ and simultaneous rapid growth of $u^{\prime}_{i}(1),$ $i=1,2.$ It is easy to see that
\begin{equation}\label{Inequality_for_Troeschs_problem}
u^{\prime}_{1}(t) = u_{2}(t) = \int\limits_{0}^{t}u^{\prime}_{2}(\xi)d\xi \leq u^{\prime}_{2}(t),\; \forall t \in [0,1].
\end{equation}
The latter means that for the Troesch's problem the 2-swap transformation on $(\varepsilon, 1)$ would be preferable over the 1-swap, in order to suppress the stiffness. The result of a 2-swap operator applied to problem \eqref{Troesch_problem} looks as follows  
\begin{equation}\label{Troesch_problem_left_part_v1}
\left[\begin{array}{l}
u_{1}^{\prime}(t) \\
u_{2}^{\prime}(t)
\end{array}\right] = \mathbf{T}_{0}(u_{1}(t), u_{2}(t), \lambda),\; [0, \varepsilon],
\end{equation}
\begin{equation}\label{Troesch_problem_right_part_v1}
\left[\begin{array}{l}
v_{1}^{\prime}(u) \\
v_{2}^{\prime}(u)
\end{array}\right] = \mathbf{T}_{1}(v_{1}(u), v_{2}(u), \lambda) \stackrel{def}{=} \frac{1}{\lambda \sinh\left(\lambda v_{1}(u)\right)}\left[\begin{array}{l}
u \\
1
\end{array}\right],\; [u_{2}(\varepsilon), u_{2}(1)],
\end{equation}
\begin{equation}\label{Matching_and_boundary_conditions_Troesch_v1}
v_{1}(u_{2}(\varepsilon)) = u_{1}(\varepsilon), \; v_{2}(u_{2}(\varepsilon))= \varepsilon,\; u_{1}(0) = 0,\; v_{1}(u_{2}(1)) = 1.
\end{equation}
It can be fairly noted that problem \eqref{Troesch_problem_left_part_v1}, \eqref{Troesch_problem_right_part_v1}, \eqref{Matching_and_boundary_conditions_Troesch_v1} looks essentially more complicated than the original one, first of all, because the reparametrized system \eqref{Troesch_problem_right_part_v1} has to be considered on an "unknown" interval $[u_{2}(\varepsilon), u_{2}(\varepsilon)]$ and, as the result, one of the boundary conditions is specified at a non-stationary (unknown) boundary point: $v_{1}(u_{2}(\varepsilon)) = 1.$  In the other words, for this particular case, the 2-swap transformed problem is, so to speak, of a different nature as compared to the initial one. With that said, we do not claim that problem \eqref{Troesch_problem_left_part_v1}, \eqref{Troesch_problem_right_part_v1}, \eqref{Matching_and_boundary_conditions_Troesch_v1} can't be efficiently tackled or that it is "bad" in the computational sense (after all the theory above clearly states that the reparametrized problem is not stiff). At this point, however, we prefer to leave this specific type of boundary value problems (with a non-stationary boundary point) for the future studies and take an alternative path.

Although, for the case of problem \eqref{Troesch_problem}, the 1-swap transformation is sub-optimal (see inequality \eqref{Inequality_for_Troeschs_problem}), its result looks more "attractive". Namely, instead of system \eqref{Troesch_problem_right_part_v1} and boundary/matching conditions \eqref{Matching_and_boundary_conditions_Troesch_v1} we get:
\begin{equation}\label{Troesch_problem_right_part_v2}
\left[\begin{array}{l}
v_{1}^{\prime}(u) \\
v_{2}^{\prime}(u)
\end{array}\right] = \frac{1}{v_{2}(u)}\left[\begin{array}{l}
1 \\
\lambda \sinh\left(\lambda u\right)
\end{array}\right],\; [u_{1}(\varepsilon), 1],
\end{equation}
\begin{equation}\label{Matching_and_boundary_conditions_Troesch_v2}
v_{2}(u_{1}(\varepsilon)) = u_{2}(\varepsilon), \; v_{1}(u_{1}(\varepsilon))= \varepsilon,\; u_{1}(0) = 0,\; v_{1}(1) = 1.
\end{equation}
Problem \eqref{Troesch_problem_left_part_v1}, \eqref{Troesch_problem_right_part_v2}, \eqref{Matching_and_boundary_conditions_Troesch_v2}, in contrast to problem \eqref{Troesch_problem_left_part_v1}, \eqref{Troesch_problem_right_part_v1}, \eqref{Matching_and_boundary_conditions_Troesch_v1}, has both its boundary points constant. However, as it can be easily verified using, for example, numerical data from \cite{Gen_sol_of_TP}, $v_{2}^{\prime}(1)  \rightarrow \infty$ as $\lambda \rightarrow \infty,$ which means that problem \eqref{Troesch_problem_left_part_v1}, \eqref{Troesch_problem_right_part_v2}, \eqref{Matching_and_boundary_conditions_Troesch_v2} still has a boundary layer at $u=1.$ To tackle this issue, we need to introduce yet another useful operator.

For any $l\in\overline{1, n}$ we define an {\it $l$-flip operator} $\mathcal{FP}_{l}: \mathbb{S}_{n} \rightarrow \mathbb{S}_{n}$ as follows
\begin{equation}\label{system_transformed_flip}
H \stackrel{def}{=} \mathcal{FP}_{l}(F) : \mathbf{w}^{\prime}(t) = \mathbf{H}(\mathbf{w}(t), t),
\end{equation}
where 
$$\mathbf{w}(t) = \left[w_{1}(t),\ldots, w_{n}(t)\right]^{T},$$
$$\mathbf{H}(\mathbf{w}(t), t) = \left[H_{1}(\mathbf{w}(t), t), \ldots, H_{n}(\mathbf{w}(t), t)\right]^{T},$$
$$H_{i}(\mathbf{w}(t), t) = F_{i}\left(w_{1}(t),\ldots, w_{l-1}(t),\frac{1}{w_{l}(t)}, w_{l+1}(t), \ldots, w_{n}(t), t\right), \; i\neq l,$$
$$
H_{l}(\mathbf{w}(t), t) = -F_{l}\left(w_{1}(t),\ldots, w_{l-1}(t),\frac{1}{w_{l}(t)}, w_{l+1}(t), \ldots, w_{n}(t), t\right)\left(w_{l}(t)\right)^{2}.$$

As one can tell, the $l$-flip transformation consists in "flipping" the $l$-th component of the unknown vector function: $$u_{l}(t) = \frac{u_{l}(t)}{1} \rightarrow \frac{1}{w_{l}(t)},$$ hence the suggested name for the operator $\mathcal{FP}_{l}$.

Applying the $2$-flip transformation to system \eqref{Troesch_problem_right_part_v2}, we get
\begin{equation}\label{Troesch_problem_right_part_v3}
\left[\begin{array}{l}
w_{1}^{\prime}(u) \\
w_{2}^{\prime}(u)
\end{array}\right] = \mathbf{T}_{2}(w_{1}(u), w_{2}(u), \lambda) \stackrel{def}{=} \left[\begin{array}{l}
w_{2}(u) \\
\lambda \sinh\left(\lambda u\right)\left(w_{2}(u)\right)^{3}
\end{array}\right],\; [u_{1}(\varepsilon), 1],
\end{equation}
and this time the boundary and matching conditions take the following form
\begin{equation}\label{Matching_and_boundary_conditions_Troesch_v3}
w_{2}(u_{1}(\varepsilon)) = \frac{1}{u_{2}(\varepsilon)}, \; w_{1}(u_{1}(\varepsilon))= \varepsilon,\; u_{1}(0) = 0,\; w_{1}(1) = 1.
\end{equation}
It is easy to ensure (again, using the numerical data from \cite{Gen_sol_of_TP}), that both $w_{1}(1)$ and $w_{2}(1)$ tend to $0$ as $\lambda$ tends to $+\infty,$ which means that BVP \eqref{Troesch_problem_left_part_v1}, \eqref{Troesch_problem_right_part_v3}, \eqref{Matching_and_boundary_conditions_Troesch_v3} does not have boundary layers. Furthermore, it is easy to verify that applying operators $\mathcal{SP}_{1}$ and $\mathcal{FP}_{2}$ (in any order, since the operators commutate) to equation \eqref{Second_order_BVP} (written in a "system" form) results in a systems of ODEs which is equivalent to equation \eqref{Inverse_second_order_BVP}. The latter means that the approach from \cite{Makarov_Dragunov_2019} is actually a partial case of the described transformation-based methodology. 

\section{Implementation aspect}\label{Implementation_aspect_section}
This section is aimed to illustrate (in a reasonable depth) how the introduced "swap" and "flip" transformations can be integrated into a pretty much any existing numerical method for solving boundary value problems. To conform to the practical orientation of the paper, we are going to pick a concrete method as a "base" one and demonstrate how it can be modified by means of the transformation-based approach. In the next section we will compare the "stiffness resistance" capabilities of the modified method with those of the base one.

Let's consider system \eqref{system} on some interval 
$[a, b]$ and supplement it by some boundary conditions 
\begin{equation}\label{boundary_conditions_general}
\mathbf{g}(\mathbf{u}(a), \mathbf{u}(b)) = 0.
\end{equation}
When it comes to the numerical methods for solving BVPs of type \eqref{system}, \eqref{boundary_conditions_general}, it is difficult to imagine anything more simple and popular than the {\it trapezoidal scheme} (see, for example, \cite[Sect. 5.1.3 Simple schemes for nonlinear problems]{Ascher_1988}), which from now on is our "base" method. Recall that the trapezoidal scheme approximates the solution of the target BVP   \eqref{system}, \eqref{boundary_conditions_general} on a mesh
\begin{equation}\label{mesh}
a = t_{0} < t_{2} < \ldots < t_{m} = b
\end{equation}
via the solution to the system of nonlinear equations
\begin{equation}\label{trapezoidal_scheme}
\frac{\mathbf{u}_{i+1} - \mathbf{u}_{i}}{h_{i}} = \frac{1}{2}\left(\mathbf{F}(\mathbf{u}_{i+1}, t_{i+1}) + \mathbf{F}(\mathbf{u}_{i}, t_{i})\right),\; h_{i} = t_{i+1} - t_{i},\; 0\leq i < m,
\end{equation}      
\begin{equation}\label{trapezoidal_scheme_bc}
\mathbf{g}(\mathbf{u}_{0}, \mathbf{u}_{m}) = 0.
\end{equation}

As a rule, to solve system \eqref{trapezoidal_scheme}, \eqref{trapezoidal_scheme_bc} some iterative procedure (like the Newton's method) has to be used. The latter implicitly assumes availability of an initial guess 
\begin{equation}\label{initial_guess}
\mathbf{u}_{i}^{(0)}, \; 0\leq i < m,
\end{equation}
which provides a reasonably good (for the iterative procedure to converge) approximation of the unknown solution.

What follows should not be taken as a strict set of instructions on how the transformation-based approach should be implemented but rather as a guidance which, by the way, follows the open source C++ implementation available at \url{https://github.com/imathsoft/MathSoftDevelopment} \footnote{To be more specific, we mean {\bf class trapezoidal\_solver}, which can be found in the repository together with a set of unit tests associated with it.}.

There are two main questions that need to be answered in scope of this section:
\begin{itemize}
\item How, in principle, the "swap" and "flip" transformations can be incorporated into the framework of the trapezoidal scheme?
\item How to decide on what transformations should be applied in order to suppress the stiffness?
\end{itemize} 

We start with the first question. In principle, transformations, as such, can be applied to system \eqref{system} on a  "sub-interval basis". Namely, for each particular sub-interval $[t_{i}, t_{i+1}],$ we can come up with a transformation $\mathcal{T}$ (some composition of the "swaps" and "flips") which, in effect, will result in the $i$-th equation from \eqref{trapezoidal_scheme}  being substituted by 
\begin{equation}\label{trapezoidal_scheme_transformed_eq}
\frac{\mathbf{q}_{i+1} - \mathbf{q}_{i}}{\tau_{i+1} - \tau_{i}} = \frac{1}{2}\left(\mathcal{T}\left(\mathbf{F}\right)(\mathbf{q}_{i+1}, \tau_{i+1}) + \mathcal{T}\left(\mathbf{F}\right)(\mathbf{q}_{i}, \tau_{i})\right).
\end{equation}
Let's see how $\mathbf{q}_{i}$ relates to $\mathbf{u}_{i}.$ To do so we need to make some assumptions about the structure of operator $\mathcal{T}.$ Obviously, there is no reason to apply more than one "swap" transformation on the same interval, since, provided that the interval is narrow enough, there will always be an index $k\in \overline{1,m}$ such that the $k$-swap transformation is not "worse" in suppressing the stiffness than  any other $i$-swap transformation for $k\neq i$ (this directly follows from the definition of the "swap" operator \eqref{system_transformed_swap}). There is also not much sense in composing a $k$-swap with a $k$-flip transformation, since, while both of them can "suppress" derivative of $u_{k}(t),$ the former have more general "impact" and thus is more preferable. Finally, it is obvious that
$$\mathcal{SP}_{k}\circ \mathcal{FP}_{l} = \mathcal{FP}_{l}\circ \mathcal{SP}_{k},\; k\neq l,$$
$$\mathcal{FP}_{l}\circ \mathcal{FP}_{k} = \mathcal{FP}_{k}\circ \mathcal{FP}_{l},\; \mathcal{FP}_{k}\circ \mathcal{FP}_{k} = \mathcal{I},$$
where sign $\circ$ denotes the composition operation and $\mathcal{I}$ is the identity operator. In the light of the above, without loss of generality, we may assume that 
\begin{equation}\label{operator_T}
\mathcal{T} = \mathcal{FP}_{1}\circ \mathcal{FP}_{2}\circ \ldots \circ \mathcal{FP}_{k-1} \circ \mathcal{SP}_{k},\; k \in \overline{1, m},
\end{equation}
in which case, components of the unknown vector $\mathbf{q}$ can be expressed through the components of $\mathbf{u}_{i}$ as follows
\begin{equation}\label{q_relates_u}
\mathbf{q}_{i} = \left[\frac{1}{u_{i,1}}, \frac{1}{u_{i, 2}}, \ldots, \frac{1}{u_{i, k-1}}, t_{i}, u_{i,k+1}, \ldots, u_{i,n}\right]^{T},\; \mathbf{u}_{i} = \left[u_{i,1}, u_{i,2}, \ldots, u_{i, n}\right]^{T}.
\end{equation}
and 
\begin{equation}\label{tau_relates_u}
\tau_{i} = u_{i, k}.
\end{equation}
Apparently, a transformation (when applied in the way described above) does not alter the number of nonlinear equations and thus, for the approach to work, it must not change the number of unknowns which for system \eqref{trapezoidal_scheme} is equal to $m\times n.$ Equalities \eqref{q_relates_u}, \eqref{tau_relates_u} demonstrate that the number of unknowns associated with the $i$-th point of the mesh indeed does not change, provided that we agree upon whether it $t_{i}$ or $\tau_{i}$ should be treated as an "unknown". Transformed equation \eqref{trapezoidal_scheme_transformed_eq} "operates" within its "natural" set of unknowns $\mathbf{q}_{i}, \mathbf{q}_{i+1}$ and from its perspective, $\tau_{i}, \tau_{i+1}$ are "fixed" points on the mesh. At the same time an equation associated with the interval to the left from the $i$-th mesh point has its own set of "natural" unknowns assuming $t_{i}$ to be a "fixed" point defined by the mesh. To solve this dilemma we can adopt a convention stating that the set of "actual" unknowns associated with the $i$-th point of the mesh coincide with the set of "natural" unknowns of the equation associated with the mesh interval 
\begin{itemize}
\item preceding the $i$-th point, if $i>0$;
\item following the $i$-th point, if $i=0$.
\end{itemize}
Although this is not the only possible convention we can embrace, it is the one used in the C++ implementation mentioned above. According to it, $t_{i}$ is a fixed point on the mesh and $\tau_{i}$ is an unknown (provided that $i > 0$). At the same time, with respect to the point with index $(i+1)$, $\tau_{i+1}$ is "fixed" and $t_{i+1}$ is an "unknown".

Another important point to emphasize is that, in contrast to the classical trapezoidal scheme, its transforma-tion-based "extension" can alter the mesh. Furhtermore, as one can notice from the reasoning above, when transformations are involved, the very definition of the mesh as a set of fixed values that the "independent" variable can take \eqref{mesh}, does not make much sense any more (since on each sub-interval the "independent" variable can be different). Instead, the mesh should be thought of as a result of an iterative procedure applied to solve the "transformed" nonlinear system. Namely, on the iteration with index $j+1,$ the mesh $\mu^{(j+1)}$ is the combination of the classical mesh \eqref{mesh} (which is subject to change after each iteration) and the approximation $\left\{\mathbf{u}_{i}^{(j)}\right\}_{i=0}^{m}$ obtained by the previous iteration: 
\begin{equation}\label{evolving_mesh}
\mu^{(j+1)} = \left\{\left(\mathbf{u}^{(j)}_{0}, t^{(j)}_{0}\right), \left(\mathbf{u}^{(j)}_{1}, t^{(j)}_{1}\right), \ldots, \left(\mathbf{u}^{(j)}_{n}, t^{(j)}_{n}\right)\right\},\; t^{(j)}_{i} < t^{(j)}_{i+1},\; i\in \overline{0, n-1}. 
\end{equation}

Now we finally can address the second question posed above.
Obviously, the only reason why we should consider applying any transformations to system \eqref{system} (at least in the current context) is to suppress stiffness of BVP \eqref{system}, \eqref{boundary_conditions_general}. According to the definition of stiffness adopted in the present paper (see section \ref{introduction_section}), the necessary condition for a (nonidentity) transformation to be applied on interval $[t^{(j)}_{i}, t^{j}_{(i+1)}]$ can be expressed as follows
\begin{equation}\label{stiffnes_inequality}
\alpha \|\mathbf{F}(\mathbf{u}^{(j)}_{i}, t^{(j)}_{i})\| + \beta \|\mathbf{F}(\mathbf{u}^{(j)}_{i+1}, t^{(j)}_{i+1})\| \geq \Theta > 1,\; \alpha + \beta = 1, \; \alpha, \beta \geq 0,
\end{equation}
where constant $\Theta$ plays a role of the "stiffness tolerance", which can be chosen individually for each problem and for each "base" method.

First of all, when facing situation \eqref{stiffnes_inequality}, we should consider applying a $k$-swap transformation, since, as it was shown in the previous section, the latter has an ultimate ability to "neutralize" stiffness, provided that $k$ is calculated as follows
\begin{equation}\label{formula_for_index_k}
k = \argmax\limits_{i\in \overline{0, n-1}}\left( \alpha |F_{i}(\mathbf{u}^{(j)}_{i+1}, t^{(j)}_{i+1})| + \beta |F_{i}(\mathbf{u}^{(j)}_{i+1}, t^{(j)}_{i+1})|\right).
\end{equation}

Depending on the boundary conditions, it can happen that for the index $k$ chosen according to formula \eqref{formula_for_index_k}, the resulted "$k$-swap transformed" BVP will have (at least) one of its boundary points being an unknown function (this situation was illustrated in the previous section on the example of the Troesch's problem). Earlier in this paper we agreed not to consider such types of problems (i.e., those with non-stationary boundary points) and to leave them for the future studies. So, for now, we assume that the $\argmax$ is taken over the subset of indices $\mathbb{I}^{\ast}\subset \overline{1, n-2}$ (i.e., the boundary intervals are excluded from the consideration). Then it also may happen that the  $k$-swap transformation that we came up with is a "suboptimal" one, meaning that the "level of stiffness" of the transformed problem is still above the acceptable tolerance threshold $\Theta,$ i.e,
\begin{equation}\label{stiffness_persists_condition}
\alpha \left\|\mathcal{SP}_{k}\left(\mathbf{F}(\mathbf{u}^{(j)}_{i}, t^{(j)}_{i})\right)\right\| + \beta \left\|\mathcal{SP}_{k}\left(\mathbf{F}(\mathbf{u}^{(j)}_{i+1}, t^{(j)}_{i+1})\right)\right\| \geq \Theta.
\end{equation}
In this case we have to consider using $l$-flip transformations for those indexes $l\in \overline{0, n-1} \setminus \mathbb{I}^{\ast}$ that cause the "violation" \eqref{stiffness_persists_condition}. As it follows from equalities \eqref{system_transformed_flip}, the necessary condition for an $l$-flip transformation to be efficient on interval $[t_{i}, t_{i+1}]$ can be formulated, for example, in the following way
$$\min\left\{|u_{i, l}|, |u_{i+1, l}|\right\} > 1.$$

We conclude the present section by briefly touching another important aspect of the transformation-based approach which is the {\it mesh refinement}. There are two possible consequences of the fact that mesh \eqref{evolving_mesh} "evolves" at each iteration, namely
\begin{itemize}
\item the appearance of "zigzags" in the mesh, i.e., situations when $t^{(j)}_{i} < t^{(j)}_{r}$ for $i > r;$
\item an extreme mesh "condensation" around some points and, as the result, a substantial "exhaustion" of the mesh on the adjacent regions.
\end{itemize}
In practice the former issue can be successfully solved by sorting the mesh points \eqref{evolving_mesh} with respect to $t^{(j)}_{i}$ values. As for the latter problem --- it also can be pretty much easily fixed by "decimating" the dense regions (i.e., by removing mesh points if they are closer than some acceptable threshold) as well as by generating new mesh points to fill the gaps in sparse regions (via linear or nonlinear interpolation). In general, when refining the mesh (by adding/removing points) one can consider using any of the known step size selection strategies appropriate for the "base" method in hands (see, for example, \cite{step_size_strategies_2012}). The strategies, of course, should be applied with respect to the actual independent variables on each sub-interval.

In the numerical examples below we use a quite simple mesh refinement strategy that is based on the following requirements:
\begin{itemize}
\item step size on the $i$-th interval, i.e., $$h_{i} = |\tau_{i+1} - \tau_{i}|$$ should be maximized under the restriction  
\begin{equation}\label{mesh_refinement_inequality}
	\left\|\frac{\mathbf{q}_{i+1} - \mathbf{q}_{i}}{h_{i}} - \mathbf{q}^{\prime}(\tau_{i})\right\| < M;
\end{equation}
\item step size $h_{i}$ should be kept not less than some minimal value $h_{min} > 0,$ and not greater than some maximal value $h_{max} > h_{min}.$
\end{itemize}
Values $M,$ $h_{min}$ and $h_{max}$ are the input data for the refinement algorithm. In formula \eqref{mesh_refinement_inequality} $\mathbf{q}(\tau)$ denotes the exact solution of the corresponding transformed system of ODEs. By decreasing $M$ we force the refinement procedure to produce more mesh points in the regions where $\|\mathbf{q}^{\prime}(\tau)\|$ is "high". At the same time, the restriction with $h_{\min}$ ensures that we do not end up with the mesh being "too fine", which, provided that we work within a finite precision arithmetic, would result in higher approximation errors of the corresponding finite differences operators. In practice, condition \eqref{mesh_refinement_inequality} can be reduced to
$$\left\|\mathcal{T}\left(\mathbf{F}\right)(\mathbf{q}_{i+1}, \tau_{i+1}) - \mathcal{T}\left(\mathbf{F}\right)(\mathbf{q}_{i}, \tau_{i})\right\| < 2 M,$$
since
$$\frac{\mathbf{q}_{i+1} - \mathbf{q}_{i}}{h_{i}} - \mathbf{q}^{\prime}(\tau_{i}) = \frac{1}{2}\mathbf{q}^{\prime\prime}(\tau^{\ast}) h_{i} \approx \frac{1}{2}\left(\mathcal{T}\left(\mathbf{F}\right)(\mathbf{q}_{i+1}, \tau_{i+1}) - \mathcal{T}\left(\mathbf{F}\right)(\mathbf{q}_{i}, \tau_{i})\right), \; \tau^{\ast} \in [\tau_{i}, \tau_{i+1}].$$
 
\section{Numerical experiments}\label{Numerical_experiments_section}
The present section is aimed to illustrate how the transformation-based approach can be applied to concrete stiff BVPs and to compare its performance with that of the corresponding "base" method (which, as we agreed above, is the {\it trapezoidal scheme} \cite{Ascher_1988}).

\subsection{Preliminary comments}
To be able to do the comparison, we need to agree about the methodology/protocol for assessing the "stiffness resistance" property of a method in a quantitative way. As a rule, stiff problems depend on some parameter, $\lambda$ which is, in a way, "proportional" to the stiffness of the problem. With this in mind, we are going to use the following protocol for evaluating the "stiffness resistance" capability of a method:
\begin{enumerate}
\item pick a mesh refinement procedure, its parameters (including the "identity" mesh refinement if we want to work with fixed meshes) and a lambda increment step $\Delta \lambda$;
\item come up with an initial guess for the solution of the BVP in question and use it to solve the problem for $\lambda = \lambda_{0}$ (the parameter should be chosen low enough for the method to succeed/converge); this constitutes the $0$-th iteration; 
\item on the $n$-th iteration, use the solution of the problem with $\lambda = \lambda_{n-1}$ (from the iteration $n-1$) as an initial guess to solve the problem with  $\lambda = \lambda_{n} = \lambda_{n-1} + \Delta \lambda;$
\item continue the process until a certain {\it stop criteria} is met (for example, the Newton's iterations does not converge or the desired precision cannot be achieved);
\item the maximal value of $\lambda_{n}$ that we can achieve without meeting the {\it stop criteria} is, obviously, a quantitative measure of the method's ability to resist stiffness; the number $\lambda_{n}$ will be called a {\it stiffness resistance number} (SRN) of the method in the context of the given problem and the given mesh refinement procedure.   
\end{enumerate}

In what follows we will refer to the "convergence" and "accuracy" stop criteria. The former criteria is met when the underlying iteration procedure fails to converge, while the latter criteria is met when the method fails to approximate the solution with some "acceptable" accuracy. Obviously, SRN calculated in terms of "accuracy" is always lower or equal to that calculated with respect to the "convergence" stop criteria. 

\begin{figure}[h!]
\centering
\begin{subfigure}[b]{0.47\textwidth}
\begin{tikzpicture}
\begin{axis}[no markers,
    every axis plot/.append style={ultra thick}, xlabel={$t$}, ylabel={$u_{1}(t)$}, legend pos= north west]
\addplot table [x index = {2}, y index = {0}, col sep=space] {Troesch_solution_lambda_3.csv};
\addplot table [x index = {2}, y index = {0}, col sep=space] {Troesch_solution_lambda_7.csv};
\addplot table [x index = {2}, y index = {0}, col sep=space] {Troesch_solution_lambda_10.csv};
\addplot table [x index = {2}, y index = {0}, col sep=space] {Troesch_solution_lambda_20.csv};
\addplot[cyan] table [x index = {2}, y index = {0}, col sep=space] {Troesch_solution_lambda_30.csv};
\addplot table [x index = {2}, y index = {0}, col sep=space] {Troesch_solution_lambda_40.csv};
\legend { $\lambda = 3$, $\lambda = 7$, $\lambda = 10$, $\lambda = 20$, $\lambda = 30$, $\lambda = 40$}
\end{axis}
\end{tikzpicture}
\end{subfigure}
\begin{subfigure}[b]{0.47\textwidth}
\begin{tikzpicture}
\begin{axis}[no markers,
    every axis plot/.append style={ultra thick}, xlabel={$t$}, ylabel={$u_{2}(t)$}, legend pos= south east, ymode = log]
\addplot table [x index = {2}, y index = {1}, col sep=space] {Troesch_solution_lambda_3.csv};
\addplot table [x index = {2}, y index = {1}, col sep=space] {Troesch_solution_lambda_7.csv};
\addplot table [x index = {2}, y index = {1}, col sep=space] {Troesch_solution_lambda_10.csv};
\addplot table [x index = {2}, y index = {1}, col sep=space] {Troesch_solution_lambda_20.csv};
\addplot[cyan] table [x index = {2}, y index = {1}, col sep=space] {Troesch_solution_lambda_30.csv};
\addplot table [x index = {2}, y index = {1}, col sep=space] {Troesch_solution_lambda_40.csv};
\end{axis}
\end{tikzpicture}
\end{subfigure}
\caption{The Troesch's problem \eqref{Troesch_problem}. Graphs of $u_{1}(t)$ (to the left) and $u_{2}(t)$ (to the right) for different values of $\lambda$.} 
\label{Example_1_Troesch_solution_graphs}
\end{figure}
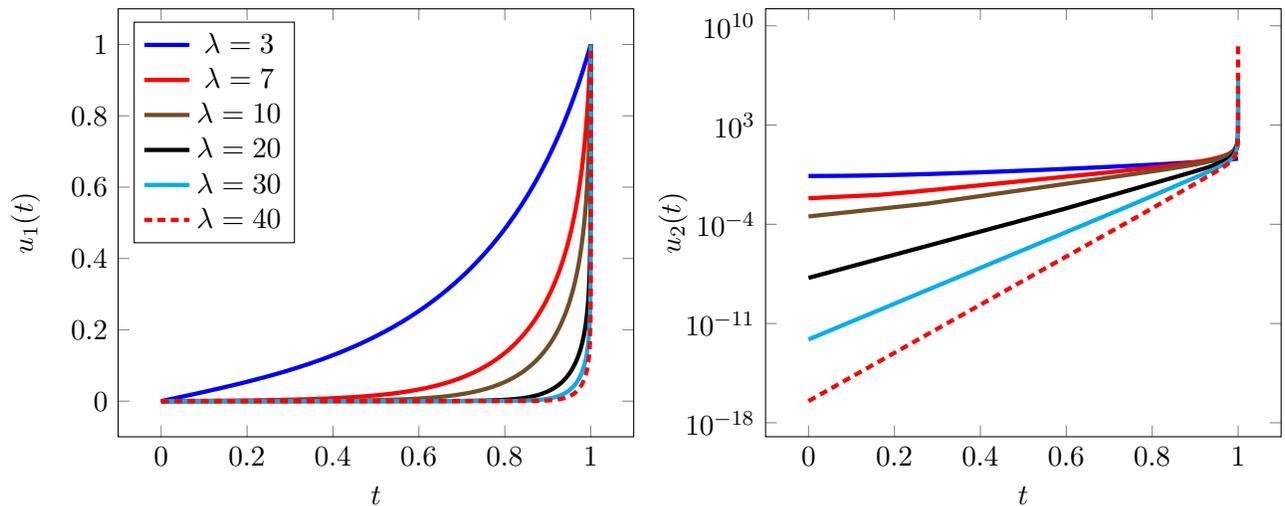

\subsection{Problem.} The two-point BVP that we will be focused on throughout this section has been already introduced above, and it is the Troesch's problem \eqref{Troesch_problem}, \cite{Troesch1976279}. Although it is definitely not the only stiff problem available (see, for example, the library of stiff BVPs collected in \cite{Cash_bvpSolve_test_problems_2010}) our choice is dictated mainly by the availability of the reference data for the Troesch's problem that we can use to validate the results of the experiments below (see \cite{Gen_sol_of_TP} and the references therein).

Figure \ref{Example_1_Troesch_solution_graphs} gives a general understanding of how a solution to the Troesch's problem \eqref{Troesch_problem} depends on the value of parameter $\lambda.$ As one can see, the higher value of $\lambda$ we choose the closer graph of $u_{1}(t)$ approaches the corner line $(0,0)$ --- $(1, 0)$  --- $(1,1).$ Consequently, for high values of $\lambda,$ function $u_{2}(t),$ which is nothing else but the derivative of $u_{1}(t),$ takes extremely low (close to zero) values almost everywhere on the interval $[0, 1]$ except for some narrow vicinity of point $t=1,$ where, in contrast, it rapidly increases, taking extremely high values. The mentioned vicinity is knows as the {\it boundary layer} of the Troesch's problem. 

\subsection{Reference data.} Before proceeding to the actual numerical experiments, it is important to clarify the origins of the reference data used in this section to assess approximation errors of the methods. A fair amount of the "ground truth" data for the Troesch's problem can be found in \cite{Gen_sol_of_TP}, namely, there one can find values of $u_{2}(0)$ and $u_{2}(1)$ calculated with up to 10 digits precision for some integer values of $\lambda$. Unfortunately, starting from $\lambda = 30$ the data provided by \cite{Gen_sol_of_TP} contains rather big "gaps": the range $30 < \lambda < 50$ (which, as it has eventually turned out, we are especially interested in) is not covered at all. To overcome this issue, the reference values of $u_{2}(0)$ and $u_{2}(1)$ for $\lambda = 1,2,\ldots 48$ have been calculated using the $\tanh$-transformation approach (see, for example, \cite{Chang20103303}, \cite{Chang20103043}) by which the Troesch's problem gets converted into a considerably less computationally difficult BVP with a polynomial type nonlinearity. The latter was solved by \verb"dsolve" subroutine within {\bf Maple 2021.1} environment (with the corresponding set of settings needed to ensure correct approximation for at least 12 significant digits in the values of $u_{2}(0)$ and $u_{2}(1)$). All the attempts to proceed past $\lambda = 48$ using the mentioned approach turned out to be unsuccessful (due to extreme time consumption). Hence, all the reference data that is used in this section for $\lambda > 48$ comes from \cite{Gen_sol_of_TP} and is limited to values of $u_{2}(0)$ for $\lambda = 50, 100, 200, 300, 400, 500.$        

\subsection{Regular trapezoidal scheme.} 
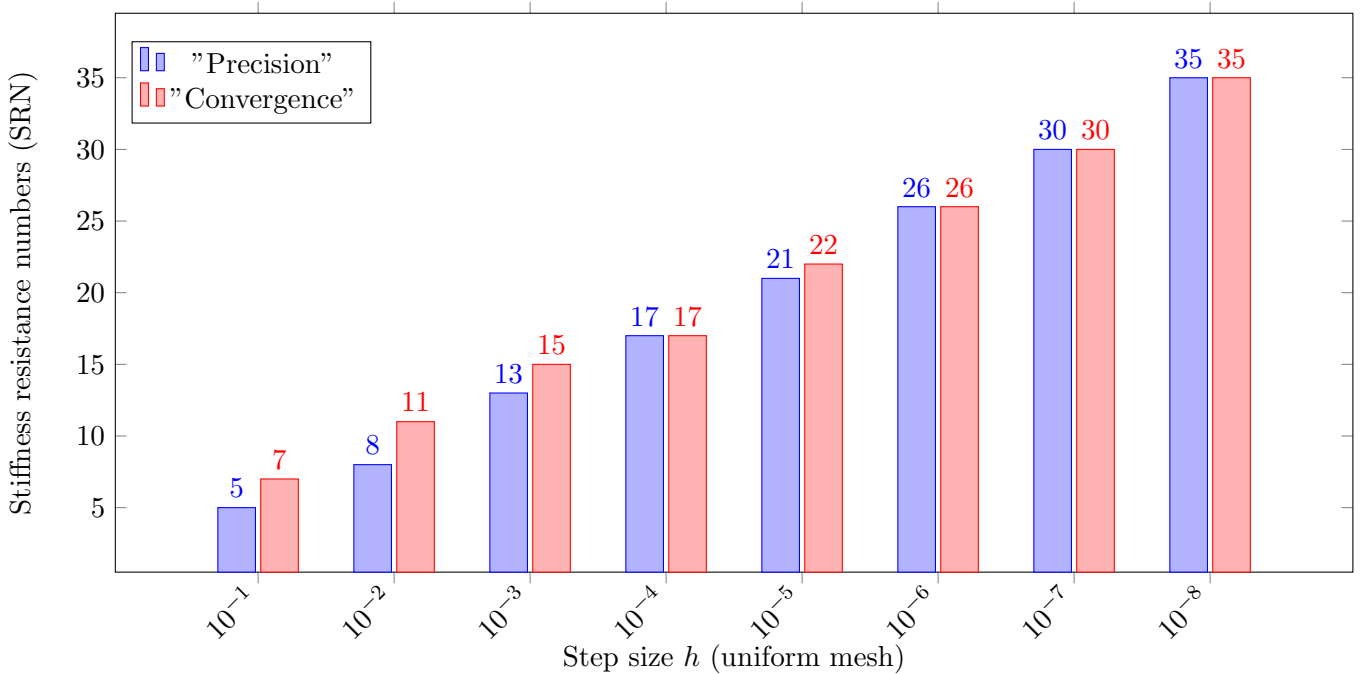
\begin{figure}
\centering
\begin{tikzpicture}
\begin{axis}[
ybar,
bar width=.5cm,
width=\textwidth,
height=.5\textwidth,
enlargelimits=0.15,
legend style={at={(0.11,0.95)}, anchor=north},
ylabel={Stiffness resistance numbers (SRN)},
xlabel = {Step size $h$ (uniform mesh)},
xlabel style={at={(0.5, -0.04)}},
symbolic x coords={$10^{-1}$, $10^{-2}$,$10^{-3}$,$10^{-4}$,$10^{-5}$,$10^{-6}$,$10^{-7}$,$10^{-8}$},
x tick label style={rotate=45, anchor=east, align=right,text width=3.5cm},
xtick=data,
nodes near coords,
nodes near coords align={vertical},
]
\addplot coordinates {($10^{-1}$,5) ($10^{-2}$,8) ($10^{-3}$,13) ($10^{-4}$,17) ($10^{-5}$,21) ($10^{-6}$,26) ($10^{-7}$,30) ($10^{-8}$,35)};
\addplot coordinates {($10^{-1}$,7) ($10^{-2}$,11) ($10^{-3}$,15) ($10^{-4}$,17) ($10^{-5}$,22) ($10^{-6}$,26) ($10^{-7}$,30) ($10^{-8}$,35)};
\legend {"Precision", "Convergence"};
\end{axis}
\end{tikzpicture}
\vspace{-65pt}
\caption{The Troesch's problem. Stiffness resistance numbers (SRN) calculated for the regular trapezoidal scheme on uniform meshes with different step sizes. Different colors correspond to different "stop criteria" used when calculating the SRNs.}\label{SRN_Troesch_uni_mesh}
\end{figure}

To begin with, we explore the stiffness resistance capabilities of the base method, the trapezoidal scheme, in its most simple version --- on uniform meshes (no mesh refinement). The corresponding SVNs calculated for different step sizes $h$, using $\lambda_{0} = 3,$ $\Delta \lambda = 1,$ are presented on Fig.~\ref{SRN_Troesch_uni_mesh}. Red bars correspond to the "convergence" stop criteria, whereas blue bars correspond to the "accuracy" stop criteria (the "acceptable" accuracy is assumed to be achieved if the method manages to approximate $u_{2}(0)$ or $u_{2}(1)$ with relative errors less than $1.0$). As one can see the SRNs calculated with the two criteria are very close to each other, so in what follows we are going to stick to the "accuracy" one, unless otherwise stated. 

\begin{figure}[h!]
\centering
\begin{subfigure}[b]{0.47\textwidth}
\begin{tikzpicture}
\begin{axis}[    
    xlabel={$\lambda$}, ylabel={$\Delta_{rel}$}, ymode = log, legend pos= south east ]
\addplot coordinates { (3, 0.000536340114205402) (6, 0.0103552948766081) (9, 0.800076971494039) };
\addplot coordinates { (3, 5.36121494049342e-06) (6, 9.96661534045059e-05) (9, 0.000971628518468911) (12, 0.0164675978295903) (15, 0.0242174906668711)};
\addplot coordinates { (3, 5.36119307237071e-08) (6, 9.96270298343716e-07) (9, 9.52233095575025e-06) (12, 7.42654566773762e-05) (15, 0.000649109521515586) (17, 0.031047960053147) (18, 1.97336038099024)};
\addplot coordinates { (3, 5.36119420168693e-10) (6, 9.96266410835925e-09) (9, 9.52038728788403e-08) (12, 7.37033145860839e-07) (15, 5.12015580578269e-06) (18, 3.52243810564472e-05) (21, 0.00203026612311887) (22, 0.922334799138397)};
\addplot coordinates { (3, 5.36121157578261e-12) (6, 9.96266159580077e-11) (9, 9.52036708789505e-10) (12, 7.36976850248746e-09) (15, 5.10774623976608e-08) (18, 3.28814148482035e-07) (21, 2.04157823853e-06) (24, 1.96189382983185e-05) (26, 0.0174866988856744)};
\addplot coordinates { (3, 5.36425204252736e-14) (6, 9.96519348165433e-13) (9, 9.52071853850922e-12) (12, 7.36980992708476e-11) (15, 5.10762650539278e-10) (18, 3.28582333097049e-09) (21, 2.00294441527163e-08) (24, 1.17793120878276e-07) (27, 7.54274821781832e-07) (30, 0.000248942265361437)};
\addplot coordinates { (3, 6.51528588161218e-16) (6, 1.00502339225373e-14) (9, 9.54418600763387e-14) (12, 7.37330123093706e-13) (15, 5.10781188490288e-12) (18, 3.2858357173123e-11) (21, 2.00256140688013e-10) (24, 1.17199580265638e-09) (27, 6.65560375965584e-09) (30, 3.7963828126766e-08) (33, 4.51954175447195e-07) (35, 0.00938496512098052)};
\end{axis}
\end{tikzpicture}
\end{subfigure}
\begin{subfigure}[b]{0.47\textwidth}
\begin{tikzpicture}
\begin{axis}[
    xlabel={$\lambda$}, ylabel={$\Delta_{rel}$}, ymode = log, legend pos= south east ]
\addplot coordinates { (3, 0.000339060397688961) (6, 0.0288665338585863) (9, 1.03429036105362)};
\addplot coordinates { (3, 3.39330724837134e-06) (6, 0.000300917286763954) (9, 0.0133887745254616) (12, 0.393870482506753) (15, 2.65538982923777)};
\addplot coordinates { (3, 3.39333435950694e-08) (6, 3.0107204787467e-06) (9, 0.000136674039548533) (12, 0.00484391549009433) (15, 0.134459412964028) (18, 2.63761599376254)};
\addplot coordinates { (3, 3.3933346926085e-10) (6, 3.0107359953704e-08) (9, 1.36705475397624e-06) (12, 4.88217865343719e-05) (15, 0.00152841994639838) (18, 0.041893399817653) (21, 0.931703764965427)};
\addplot coordinates { (3, 3.39326390852795e-12) (6, 3.01073603083413e-10) (9, 1.36705789641754e-08) (12, 4.88257969064023e-07) (15, 1.5323110443106e-05) (18, 0.000442872345281457) (21, 0.0118927135416834) (24, 0.263947810732751) (26, 1.87594303354642)};
\addplot coordinates { (3, 3.39347209700016e-14) (6, 3.01069347435649e-12) (9, 1.36705748590921e-10) (12, 4.88258373733458e-09) (15, 1.5323505559479e-07) (18, 4.4320143902353e-06) (21, 0.000121141268124379) (24, 0.00316204528339988) (27, 0.0740358077676506) (30, 1.45419970288185)};
\addplot coordinates { (3, 2.08188472208599e-16) (6, 2.99668530046674e-14) (9, 1.36699275190271e-12) (12, 4.88259317774271e-11) (15, 1.53235091561253e-09) (18, 4.43204745361188e-08) (21, 1.21165945835481e-06) (24, 3.17852294019388e-05) (27, 0.000806947287009142) (30, 0.0194602474031981) (33, 0.392908986865491) (35, 2.48467393936476)};
\legend {$h = 10^{-2}$, $h = 10^{-3}$, $h = 10^{-4}$, $h=10^{-5}$, $h=10^{-6}$, $h=10^{-7}$, $h=10^{-8}$}
\end{axis}
\end{tikzpicture}
\end{subfigure}
\caption{The Troesch's problem. Relative errors of the approximations for $u_{2}(0)$ (left) and $u_{2}(1)$ (right) calculated with the regular trapezoidal scheme on uniform meshes with different step sizes $h.$ } \label{Example_1_graph_delta_uniform}
\end{figure}
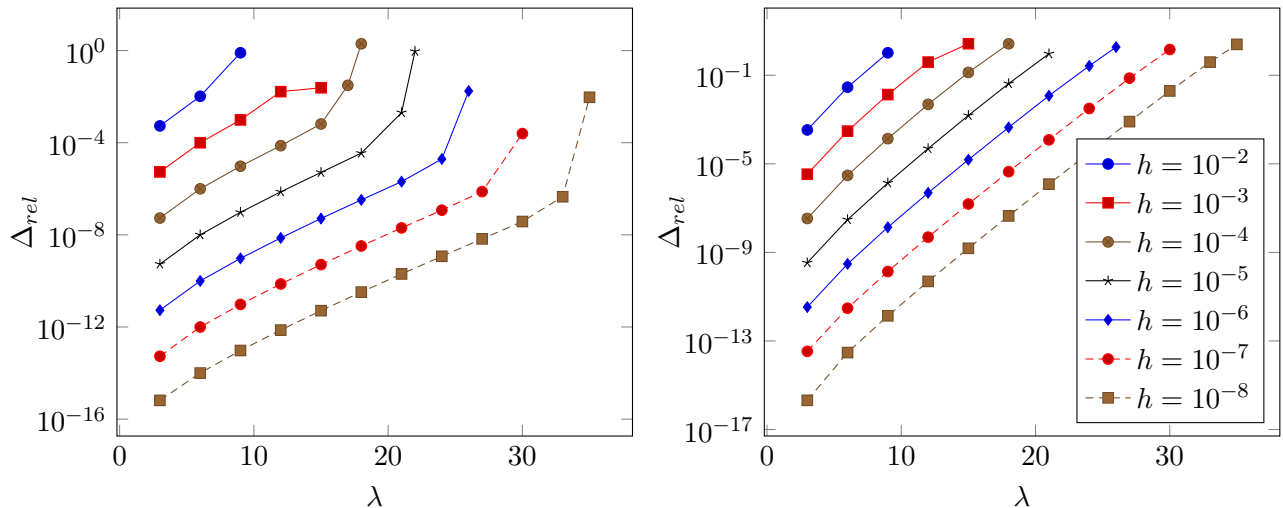

Let's asses the approximation quality of the trapezoidal scheme for different values of $\lambda$ and on different meshes. The relative errors of the approximations for $u_{2}(0),$ $u_{2}(1)$ evaluated with respect to the reference data are presented in Fig.~\ref{Example_1_graph_delta_uniform}. There are a few conclusions that we can make looking onto the error curves. First of all, the spacing between the polylines that correspond to different values of the step size $h,$ seems to be in a good agreement with the well known fact that the approximation error of the trapezoidal scheme behaves as $\mathcal{O}(h^{2}).$ Second of all, the approximation errors grow as $\lambda$ increases. Since the polylines representing relative errors are almost straight in the logarithmic scale, we can say that the growth is exponential. Finally, it is worth to point out that the slope of the polylines that correspond to the errors of $u_{1}(0)$ is lower than that of the polylines representing errors of $u_{2}(1),$ which, apparently, indicates that the approximation at $t=0$ is less sensitive to $\lambda.$ The polylines for $u_{2}(0)$ exhibit a rapid drop of accuracy near the point where the iteration procedure is close to meeting the stop criteria, which seems to be not the case for $u_{2}(1).$ This, however, might be explained by the fact that the absolute values of $u_{2}(0)$ are actually close to $0$ and, thus, even considerably low absolute approximation errors in them can result in the high relative errors (this is not the case for $u_{2}(1) \gg 1$).    

It took several hours to calculate the SRNs that corresponds to step size $h=10^{-8}$ which means that for the trapezoidal scheme on uniform meshes, parameter $\lambda = 30$ is very close to the absolute practical limit, i.e., with the available (and quite modern) computational resources we can barely proceed further. Obviously this is due to the large sizes of the uniform meshes which are inversely proportional to the corresponding step sizes.

\begin{figure}
\begin{tikzpicture}
\begin{axis}[
ybar,
bar width=.5cm,
width=0.45\textwidth,
height=.4\textwidth,
enlargelimits=0.15,
legend style={at={(0.18,0.95)}, anchor=north},
ylabel={Stiffness resistance numbers (SRN)},
xlabel = {$h_{min}$},
xlabel style={at={(0.5, -0.05)}},
symbolic x coords={$10^{-3}$, $10^{-4}$,$10^{-5}$,$10^{-6}$,$10^{-7}$,$10^{-8}$,$10^{-9}$,$10^{-10}$, $10^{-11}$},
x tick label style={rotate=45, anchor=east, align=right,text width=3.5cm},
xtick=data,
nodes near coords,
nodes near coords align={vertical},
]
\addplot coordinates {($10^{-3}$,14) ($10^{-4}$,17) ($10^{-5}$,21) ($10^{-6}$,26) ($10^{-7}$,30) ($10^{-8}$,34) ($10^{-9}$,38) ($10^{-10}$,44) ($10^{-11}$,48)};
\end{axis}
\end{tikzpicture}
\hfill
\begin{tikzpicture}
\begin{axis}[
ymode=log,
ybar,
bar width=.5cm,
width=0.45\textwidth,
height=.4\textwidth,
enlargelimits=0.15,
legend style={at={(0.18,0.95)}, anchor=north},
ylabel={Max mesh size (number of knots)},
xlabel = {$h_{min}$},
xlabel style={at={(0.5, -0.05)}},
symbolic x coords={$10^{-3}$, $10^{-4}$,$10^{-5}$,$10^{-6}$,$10^{-7}$,$10^{-8}$,$10^{-9}$,$10^{-10}$, $10^{-11}$},
x tick label style={rotate=45, anchor=east, align=right,text width=3.5cm},
xtick=data,
nodes near coords,
nodes near coords align={vertical},
]
\addplot coordinates {($10^{-3}$,336) ($10^{-4}$,1155) ($10^{-5}$,4874) ($10^{-6}$,21822) ($10^{-7}$,97402) ($10^{-8}$,428069) ($10^{-9}$,2029636) ($10^{-10}$,8819051) ($10^{-11}$,40147940)};
\end{axis}
\end{tikzpicture}
\vspace{-50pt}
\caption{The Troesch's problem. Stiffness resistance numbers (SRN) calculated for the regular trapezoidal scheme on non-uniform meshes produced by the mesh refinement algorithm with $M = 0.1,$ $h_{max} = 0.01$ and different step size thresholds $h_{min}$ (to the left). Maximal sizes of the meshes, the refinement algorithm produced in the process of calculating the corresponding SRNs (to the right).}\label{SRN_Troesch_non_uni_mesh}
\end{figure}
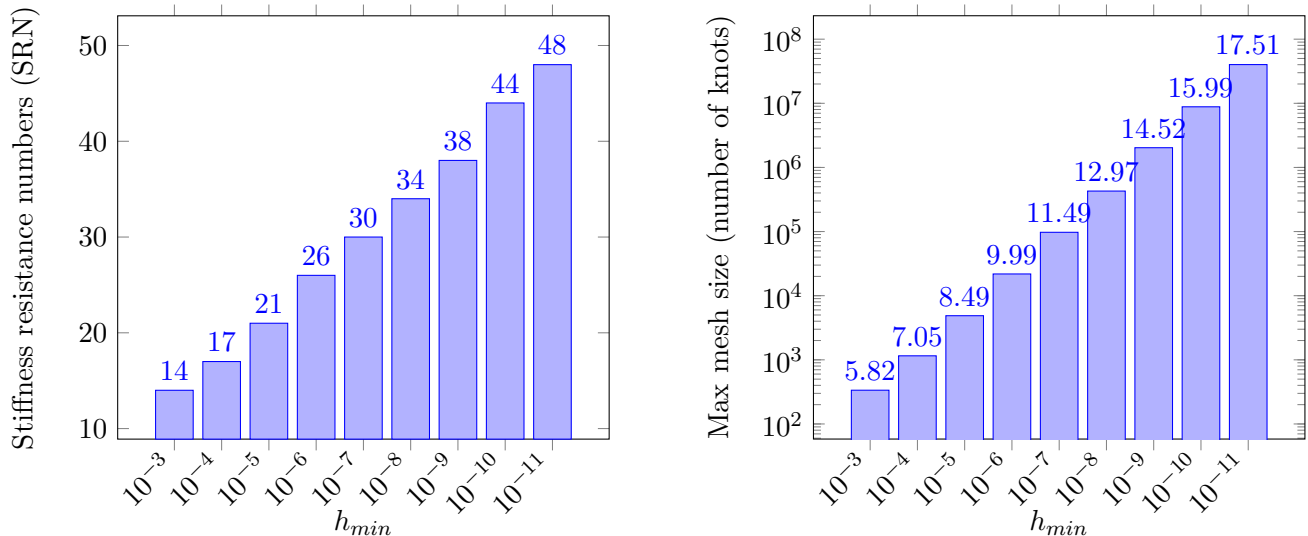

Now let's see how far we can get by using a more efficient mesh building approach, namely, the mesh refinement algorithm introduced in the previous section (see \eqref{mesh_refinement_inequality} and the related description). It is easy to see that, by design, the refinement procedure increases mesh density proportionally to the magnitude of the second derivatives of $u_{1}(t),$ $u_{2}(t)$ which means that most of the refined mesh points will be concentrated near the right boundary point $t=1,$ whereas the rest of the interval $[0,1]$ will be covered by considerably sparse mesh with distances between adjacent points determined by parameter $h_{max}.$ The results obtained with $M = 0.1$ and $h_{max} = 0.01$ are visualized in Fig.~\ref{SRN_Troesch_non_uni_mesh}. As one can see, using non-uniform meshes we managed to solve the Troesch's problem for $\lambda = 48$ (with $h_{min} = 10^{-11}$), which took a mesh containing more than 40 millions points! This also took a few hours of computations, so  $\lambda = 48$ can be considered as a practical limit for the given approach. The leftmost bar chart in Fig.~\ref{SRN_Troesch_non_uni_mesh} suggests that by decreasing $h_{min}$ 10 times we can gain about 4 more points of SRN and this will cost us about 4 times larger mesh. Following the trend, we can forecast that with $h_{min} = 10^{-16}$ we might get to $\lambda = 68,$ operating on a mesh containing about 40 billions points! The latter, however, would require us to work in the quadruple (or higher) precision arithmetic, since quantity $10^{-16}$ is already below the {\it machine epsilon} for the double-precision arithmetic (we will touch this topic in more detail a bit later in this section). 

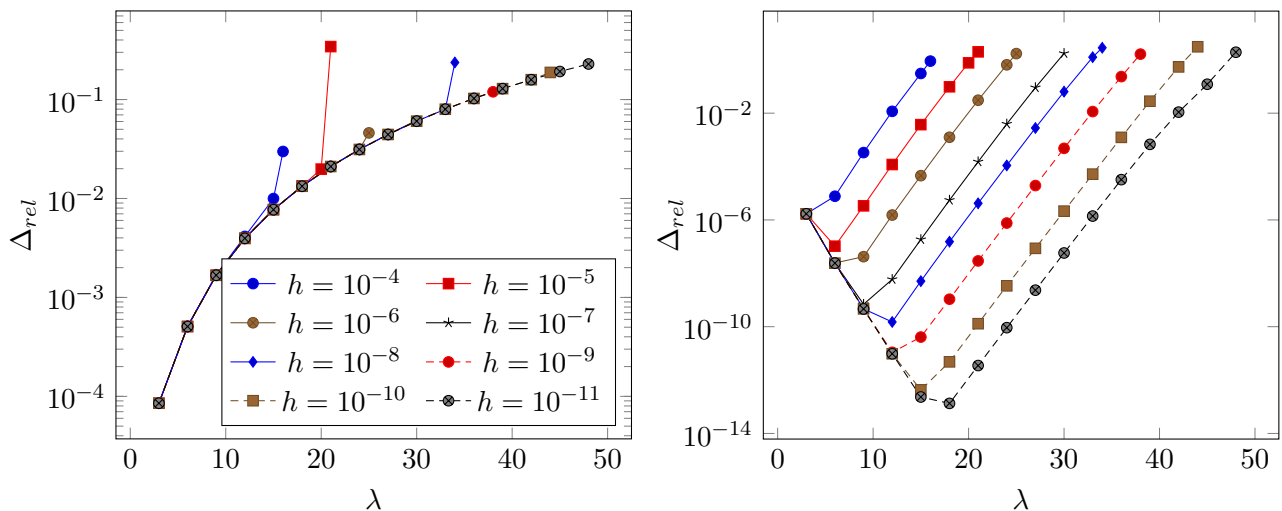
\begin{figure}[h!]
\centering
\begin{subfigure}[b]{0.47\textwidth}
\begin{tikzpicture}
\begin{axis}[ 
		legend columns=2, 
        legend style={
            /tikz/column 2/.style={
                column sep=5pt,
            },
        },   
    xlabel={$\lambda$}, ylabel={$\Delta_{rel}$}, ymode = log, legend pos= south east 
    ]
\addplot coordinates { (3, 8.52711955676491e-05) (6, 0.000510139414701257) (9, 0.00169777983884326) (12, 0.00413567348721288) (15, 0.00998482495316744) (16, 0.0298580016568288)};
\addplot coordinates { (3, 8.52711955676491e-05) (6, 0.000508163349123893) (9, 0.00167542385244996) (12, 0.00395662468705673) (15, 0.00772300765029355) (18, 0.0133904572652397) (20, 0.0198388787860109) (21, 0.342986090940559)};
\addplot coordinates { (3, 8.52711955676491e-05) (6, 0.000508149521219061) (9, 0.00167520041030954) (12, 0.00395485626446696) (15, 0.00771073386706652) (18, 0.0132977406377338) (21, 0.0210601179312343) (24, 0.0314415119934822) (25, 0.0460155609637774)};
\addplot coordinates { (3, 8.52711955676491e-05) (6, 0.000508149498039549) (9, 0.00167519843654627) (12, 0.00395483553212782) (15, 0.00771057991838642) (18, 0.0132967745286113) (21, 0.0210542885418683) (24, 0.031311335158278) (27, 0.0443618860710712) (30, 0.0615084232750576)};
\addplot coordinates { (3, 8.52711955676491e-05) (6, 0.000508149498039549) (9, 0.00167519842402351) (12, 0.00395483541941163) (15, 0.00771057927752529) (18, 0.0132967706540868) (21, 0.0210542649335824) (24, 0.0313111969639655) (27, 0.0443609821646221) (30, 0.0604832641905829) (33, 0.0798663943789401) (34, 0.236971188198129)};
\addplot coordinates { (3, 8.52711955676491e-05) (6, 0.000508149498039549) (9, 0.00167519842395277) (12, 0.00395483541772111) (15, 0.00771057926342484) (18, 0.0132967705508275) (21, 0.0210542642722829) (24, 0.031311193077462) (27, 0.0443609602567011) (30, 0.0604831317814549) (33, 0.0798468313676208) (36, 0.102625556577559) (38, 0.11976606482022)};
\addplot coordinates { (3, 8.52711955676491e-05) (6, 0.000508149498039549) (9, 0.00167519842395277) (12, 0.0039548354177092) (15, 0.00771057926329941) (18, 0.0132967705499911) (21, 0.0210542642674733) (24, 0.0313111930501534) (27, 0.0443609601045741) (30, 0.0604831309540883) (33, 0.0798468268946401) (36, 0.102625520714673) (39, 0.128906723194355) (42, 0.158741700291933) (44, 0.187989943116367)};
\addplot coordinates { (3, 8.52711955676491e-05) (6, 0.000508149498039549) (9, 0.00167519842395277) (12, 0.0039548354177092) (15, 0.00771057926329854) (18, 0.0132967705499868) (21, 0.0210542642674513) (24, 0.0313111930500334) (27, 0.0443609601039117) (30, 0.0604831309504958) (33, 0.0798468268755676) (36, 0.102625520615268) (39, 0.128906722662895) (42, 0.158741691742494) (45, 0.192064847171227) (48, 0.228745476822859)};
\legend {
$h = 10^{-4}$, $h = 10^{-5}$, $h=10^{-6}$, $h=10^{-7}$, $h=10^{-8}$, $h=10^{-9}$, $h=10^{-10}$, $h=10^{-11}$}
\end{axis}
\end{tikzpicture}
\end{subfigure}
\begin{subfigure}[b]{0.47\textwidth}
\begin{tikzpicture}
\begin{axis}[
    xlabel={$\lambda$}, ylabel={$\Delta_{rel}$}, ymode = log, legend pos= south east ]
\addplot coordinates { (3, 1.69451434856383e-06) (6, 7.72145500834107e-06) (9, 0.000333757995765411) (12, 0.0117033204438981) (15, 0.307057656147005) (16, 0.892505072792363)};
\addplot coordinates { (3, 1.69451434856383e-06) (6, 1.04000138188894e-07) (9, 3.364590348404e-06) (12, 0.000119182788161911) (15, 0.00371815077587403) (18, 0.0973134895172892) (20, 0.768447695556176) (21, 1.97914051812118)};
\addplot coordinates { (3, 1.69451434856383e-06) (6, 2.43619393501789e-08) (9, 4.19925063520067e-08) (12, 1.52407874976935e-06) (15, 4.55818023994077e-05) (18, 0.00125280577042064) (21, 0.0305795511424057) (24, 0.648729718861671) (25, 1.70433889449104)};
\addplot coordinates { (3, 1.69451434856383e-06) (6, 2.4028383983581e-08) (9, 6.90536278102177e-10) (12, 5.92383609687661e-09) (15, 1.8481993251056e-07) (18, 5.50982005188447e-06) (21, 0.000150886103918983) (24, 0.00393378681438693) (27, 0.0909569148686327) (30, 1.73690180625535)};
\addplot coordinates { (3, 1.69451434856383e-06) (6, 2.4028383983581e-08) (9, 4.66334309559343e-10) (12, 1.50049510344549e-10) (15, 5.08590086700711e-09) (18, 1.5304365539306e-07) (21, 4.1923243680632e-06) (24, 0.000109963653832792) (27, 0.00278295407328713) (30, 0.0640879878923142) (33, 1.25458540111506) (34, 2.81914286914766)};
\addplot coordinates { (3, 1.69451434856383e-06) (6, 2.4028383983581e-08) (9, 4.63453488382256e-10) (12, 1.09820126936424e-11) (15, 4.08523689595593e-11) (18, 1.07177359128519e-09) (21, 2.91149932647649e-08) (24, 7.63778216917196e-07) (27, 1.94152038206345e-05) (30, 0.000481061781886837) (33, 0.0114915111047008) (36, 0.234466911257259) (38, 1.62931746216237)};
\addplot coordinates { (3, 1.69451434856383e-06) (6, 2.4028383983581e-08) (9, 4.63453488382256e-10) (12, 9.67120497897423e-12) (15, 4.25561229603477e-13) (18, 4.81377333769955e-12) (21, 1.29412279817745e-10) (24, 3.39457531770802e-09) (27, 8.62926303283853e-08) (30, 2.13978738506355e-06) (33, 5.19999214025918e-05) (36, 0.00124049242782798) (39, 0.0281434658998917) (42, 0.540819484477403) (44, 3.07510950564336)};
\addplot coordinates { (3, 1.69451434856383e-06) (6, 2.4028383983581e-08) (9, 4.63453488382256e-10) (12, 9.66331448859248e-12) (15, 2.31769901347284e-13) (18, 1.33678512525652e-13) (21, 3.49938363584347e-12) (24, 9.17990940455244e-11) (27, 2.33296367955973e-09) (30, 5.77726863586551e-08) (33, 1.39491800577397e-06) (36, 3.23438503646756e-05) (39, 0.000678044667308139) (42, 0.0109292277943463) (45, 0.121348723400984) (48, 1.92007220994707)};
\end{axis}
\end{tikzpicture}
\end{subfigure}
\caption{The Troesch's problem. Relative error of the approximations for $u_{2}(0)$ (left) and $u_{2}(1)$ (right) calculated with the regular trapezoidal scheme on non-uniform (refined) meshes with $M = 0.1,$ $h_{max} = 0.01$ and different step size thresholds $h = h_{min}.$ } 
\label{Example_1_graph_delta_nonuniform}
\end{figure}

Our experiment with non-uniform meshes cannot be considered complete until we examine the approximation errors, just as we did it for the case of uniform meshes. The corresponding relative errors for $u_{2}(0)$ and $u_{2}(1)$ are shown in Fig.~\ref{Example_1_graph_delta_nonuniform}. As one can see, the approximation errors at the left boundary point $t = 0$ does not change much as $h_{min}$ decreases. This can be explained by the fact that, due to the specifics of the mesh refinement procedure (that have been discussed earlier in this section), the mesh near the left boundary point remains uniform and rather sparse (with the step size equal to $h_{max} = 0.01$). At the same time the corresponding error curves for $u_{2}(1)$ do demonstrate a decline as $h_{min}$ decreases, similar to how it was in the uniform meshes case. The V-shape of the curves reflects the impact of the refinement procedure on the approximation process: for lower values of $\lambda,$ adding more points within the boundary layer allows to decrease the approximation error of $u_{2}(1);$ this persists until the refinement algorithm hits the limitation of $h_{min}$ and the error starts to grow as $\lambda$ increases (exactly as it was with uniform meshes).  

\subsection{Transformation-based approach.} Now that we have some understanding about the "stiffness-resista-nce" capabilities of the regular trapezoidal scheme (with respect to the Troesch's problem) let's see how those can be improved by involving some transformations. In what follows, we will examine two transformation strategies that seem to be applicable for the problem at hands.

\subsubsection{Strategy $I$-$\mathcal{SP}_{1}\mathcal{FP}_{2}$} We start with a strategy that can be abbreviated as $I$-$\mathcal{SP}_{1}\mathcal{FP}_{2}.$ It consists in dividing the interval $[0, 1]$ onto two adjacent sub-intervals $[0, \varepsilon],$ $[\varepsilon, 1]$ with parameter $0< \varepsilon < 1$ to be deduced from the restrictions below. On the left-hand side interval $[0, \varepsilon],$ where the solution is "calm", we do not apply any transformation (hence $"I"$ in the abbreviation, standing for "identity transformation"), whereas on the right-hand side interval $[\varepsilon, 1],$ containing the boundary layer, we apply transformations $\mathcal{SP}_{1}$ and $\mathcal{FP}_{2}$ (the order is not important since they commutate with each other). All in all, the strategy results in BVP \eqref{Troesch_problem_left_part_v1}, \eqref{Troesch_problem_right_part_v3}, \eqref{Matching_and_boundary_conditions_Troesch_v3} and is totally equivalent to the transformation proposed in \cite{Makarov_Dragunov_2019}. The "matching" point $\varepsilon$ can be determined from the condition
\begin{equation}\label{first_matching_point_condition}
\varepsilon = \min\limits_{t\in [0, 1]} \left\{t : u_{2}(t) > 1\right\}.
\end{equation}
Obviously such a choice of $\varepsilon$ ensures that transformation $\mathcal{SP}_{1}\mathcal{FP}_{2}$ has a "stiffness-suppressing" effect (see the definitions of the "flip" and "swap" transformations in section \ref{transformations_section}).          

\begin{figure}
\begin{tikzpicture}
\begin{axis}[
ybar,
bar width=.3cm,
width=0.45\textwidth,
height=.4\textwidth,
enlargelimits=0.15,
legend style={at={(0.25,0.95)}, anchor=north},
ylabel={Stiffness resistance numbers (SRN)},
xlabel = {$h_{min}$},
xlabel style={at={(0.5, -0.05)}},
symbolic x coords={$0.25$, $0.1$, $0.05$, $0.025$,  $0.01$, $0.005$},
x tick label style={rotate=45, anchor=east, align=right,text width=3.5cm},
xtick=data,
nodes near coords,
nodes near coords align={vertical},
]
\addplot coordinates {($0.25$, 6) ($0.1$, 9) ($0.05$, 16) ($0.025$, 35) ($0.01$, 46) ($0.005$, 46)};
\addplot coordinates {($0.01$, 99) ($0.005$, 167)};
\legend {$h_{max} = 0.1$, $h_{max} = 0.01$};
\end{axis}
\end{tikzpicture}
\hfill
\begin{tikzpicture}
\begin{axis}[
ybar,
bar width=.3cm,
width=0.45\textwidth,
height=.4\textwidth,
enlargelimits=0.15,
legend style={at={(0.25,0.95)}, anchor=north},
ylabel={Max mesh size (number of knots)},
xlabel = {$h_{min}$},
xlabel style={at={(0.5, -0.05)}},
symbolic x coords={$0.25$, $0.1$, $0.05$, $0.025$,  $0.01$, $0.005$},
x tick label style={rotate=45, anchor=east, align=right,text width=3.5cm},
xtick=data,
nodes near coords,
nodes near coords align={vertical},
]
\addplot coordinates {($0.25$, 6) ($0.1$, 15) ($0.05$, 28) ($0.025$, 57) ($0.01$, 119) ($0.005$, 189)};
\addplot coordinates {($0.01$, 147) ($0.005$, 285)};
\legend {$h_{max} = 0.1$, $h_{max} = 0.01$};
\end{axis}
\end{tikzpicture}
\vspace{-50pt}
\caption{The Troesch's problem. Stiffness resistance numbers (SRN) calculated for the regular trapezoidal scheme enhanced with $I$-$\mathcal{SP}_{1}\mathcal{FP}_{2}$ transformation strategy and the mesh refinement algorithm with $M = 0.1$ and different step size thresholds $h_{min},$ $h_{max}$ (to the left). Maximal sizes of the meshes produced by the refinement algorithm in the process of calculating the corresponding SRNs (to the right). }\label{SRN_Troesch_transformation_old}
\end{figure}
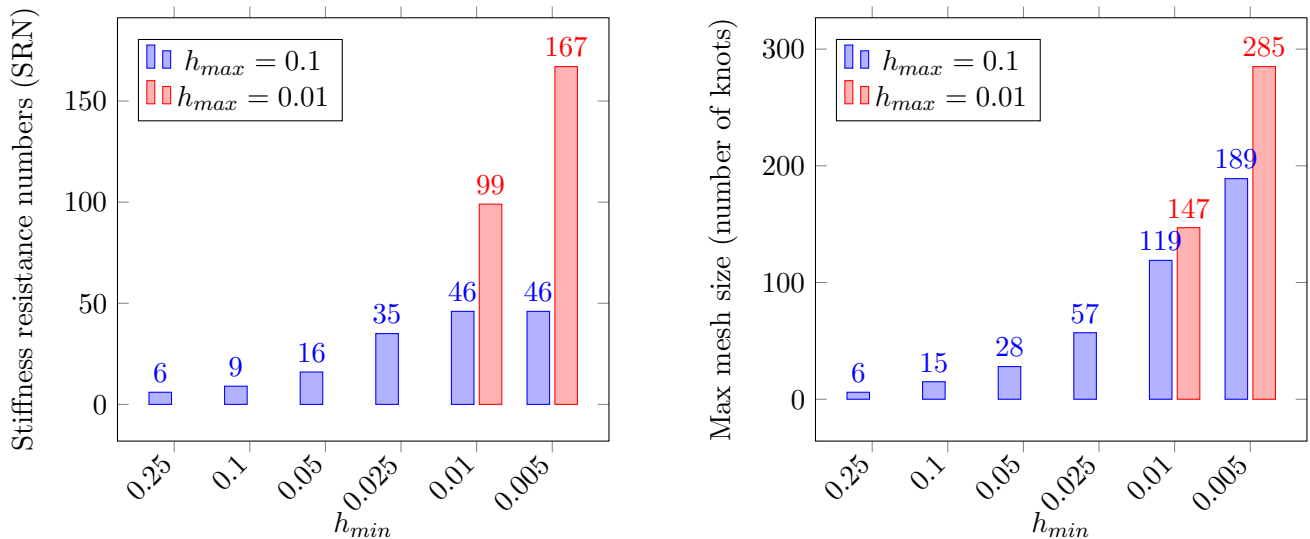

In Fig.~\ref{SRN_Troesch_transformation_old} one can find the stiffness resistance numbers calculated for $I$-$\mathcal{SP}_{1}\mathcal{FP}_{2}$ as well as the information about sizes of the meshes it took to calculate them. Comparing this data with the similar data presented in Fig.~\ref{SRN_Troesch_uni_mesh} and \ref{SRN_Troesch_non_uni_mesh}, it becomes obvious that, in terms of the introduced SRN score, $I$-$\mathcal{SP}_{1}\mathcal{FP}_{2}$ significantly outperforms the regular trapezoidal scheme on both uniform and non-uniform meshes. For example, as we can see in Fig.~\ref{SRN_Troesch_transformation_old}, approach $I$-$\mathcal{SP}_{1}\mathcal{FP}_{2}$ is able to solve the Troesch's problem for $\lambda = 46$ using a mesh containing merely $119$ points! At the same time, as we saw above, $\lambda = 46$ is practically "unreachable" for the regular trapezoidal scheme on uniform meshes and it took almost 9 millions points to get to $\lambda = 44$ using a non-uniform mesh! It is worth mentioning, that the "blue" bars in Fig.~\ref{SRN_Troesch_transformation_old} correspond to the "accuracy" stop criteria whereas "red" ones correspond to the "convergence" stop criteria (see the definition earlier in this section). The latter stop criteria was used because of the lack of reliable reference data for $\lambda > 48$ (except for the very few reference values of $u_{2}(0)$ from \cite{Gen_sol_of_TP}).

The most natural question to ask, when looking at Fig.~\ref{SRN_Troesch_transformation_old}, is probably about the "limit", i.e., the maximal value of $\lambda$ that $I$-$\mathcal{SP}_{1}\mathcal{FP}_{2}$ can handle. Besides a number of other factors, the answer essentially depends on the precision of the arithmetic that was used to implement the corresponding solving procedure. All the data in this section is obtained using the {\it double-precision} arithmetic, unless otherwise stated. For this precision, the value $\lambda = 167$ (see Fig.~\ref{SRN_Troesch_transformation_old}) is way beyond the reasonable limit and here is why. It is well known that quantity $u_{2}(1)$ is a monotonically increasing function of $\lambda.$ Somewhere near the point $\lambda = 73,$ value of $u_{2}(1)$ exceeds threshold $\varepsilon_{d}^{-1},$ where $\varepsilon_{d}$ is known as the {\it machine epsilon} for the double-precision arithmetic (the maximal positive number that can be stored in a double-precision variable, satisfying the equality $1.0 + \varepsilon_{d} = 1.0$ in the double-precision arithmetic; $\varepsilon_{d} \approx 2.22\times 10^{-16}$), which leads to a situation when for sufficiently small $\delta > 0$ $$w_{2}(1) - w_{2}(1-\delta) = 0,$$ where $w_{1}(u_{1}(t)) = t,$ see \eqref{Troesch_problem_right_part_v3}. In the other words, for $\lambda \geq 73,$ the graph of the solution $u_{1}(t)$ to the Troesch's problem \eqref{Troesch_problem} becomes so steep near the point $t=1,$ that, in the double-precision arithmetic, it can't be distinguished from a vertical line. Apparently, in this situation, the finite difference operator (that we use within the trapezoidal scheme) is not able to approximate the corresponding differential operator and, although the iterative process remains convergent for $\lambda$ as hight as $167$ (and even higher, under certain conditions), the approximation of the solution that we get becomes really poor, especially near the right boundary point $t=1.$ This can be mitigated by using, for example, a multi-precision arithmetic and we will see the corresponding results later in this section. Before that, however, in order to get a better understanding about the practical characteristics of $I$-$\mathcal{SP}_{1}\mathcal{FP}_{2}$ let's assess its approximation capabilities for $\lambda \leq 48$ (in the double-precision arithmetic).   

\begin{figure}[h!]
\centering
\begin{subfigure}[b]{0.47\textwidth}
\begin{tikzpicture}
\begin{axis}[ 
		legend columns=2, 
        legend style={
            /tikz/column 2/.style={
                column sep=5pt,
            },
        },   
    xlabel={$\lambda$}, ylabel={$\Delta_{rel}$}, ymode = log, legend pos= south east 
    ]
\addplot coordinates { (2, 0.0108593925844998) (4, 0.0907012396704175) (6, 0.277377639879787) (8, 0.647052599832164) (9, 0.82733044336328)};
\addplot coordinates {(2, 0.00430509472338007) (4, 0.0310598357790495) (6, 0.10481334155203) (9, 0.32794447524258) (12, 0.656880217617501) (15, 0.910149925209691) (16, 0.957720228687393)};
\addplot coordinates { (2, 0.000399969689754496) (4, 0.00457189001085539) (6, 0.0158174154551311) (9, 0.055547698924466) (12, 0.128397629206847) (15, 0.24187961116382) (18, 0.392916977613911) (21, 0.551881571657056) (24, 0.708972492153557) (27, 0.844293745008411) (30, 0.92756055672838) (33, 0.973688959909348) (35, 0.988973459613834)};
\addplot coordinates { (2, 0.000240144188884976) (4, 0.00127664222617801) (6, 0.00493991105280866) (9, 0.0175775736650793) (12, 0.0422890108097879) (15, 0.0820693001224397) (18, 0.139317900731248) (21, 0.216029173981859) (24, 0.308453312751689) (27, 0.414864361074531) (30, 0.529219596189187) (33, 0.641992090227836) (36, 0.746881757201407) (39, 0.83680764117003) (42, 0.905341921922006) (45, 0.952499726679541) (48, 0.980663492533534)};
\addplot coordinates { (2, 6.1896997316782e-05) (4, 0.000696768104650927) (6, 0.00242384556484987) (9, 0.00835625369078519) (12, 0.02015180465738) (15, 0.0393635070568644) (18, 0.0682188144511865) (21, 0.107159020994097) (24, 0.155975076021183) (27, 0.216839811556007) (30, 0.288009803053877) (33, 0.364937496387452) (36, 0.447673537285184) (39, 0.535433592384743) (42, 0.618489795035836) (45, 0.697020506717757) (48, 0.770001162959733)};
\addplot coordinates { (2, 3.57295323401636e-07) (4, 6.93535023791836e-06) (6, 2.8518782536201e-05) (9, 0.000103735521946946) (12, 0.000242406923559847) (15, 0.000474703034944456) (18, 0.000825244596891765) (21, 0.00131644965552315) (24, 0.0019696780918833) (27, 0.00281126148829658) (30, 0.00386319999775533) (33, 0.00515124121316324) (36, 0.00669420800013967) (39, 0.008513565092951) (42, 0.0106438403260382) (45, 0.0130896473238096) (48, 0.0158754711496254)};
\addplot coordinates { (2, 1.39474274957041e-08) (4, 1.35845251336037e-07) (6, 4.9257748461729e-07) (9, 1.76121295999208e-06) (12, 4.3312220242653e-06) (15, 8.64966737105203e-06) (18, 1.51757822318432e-05) (21, 2.43089330717763e-05) (24, 3.65840606079041e-05) (27, 5.23077396491831e-05) (30, 7.17431805149499e-05) (33, 9.5363447627174e-05) (36, 0.000123627734465342) (39, 0.000157252305697288) (42, 0.000196378358042881) (45, 0.000241522909876783) (48, 0.000293192857207176)};
\addplot coordinates { (2, 9.87878243577411e-11) (4, 1.1118937175789e-09) (6, 4.15023192206685e-09) (9, 1.51810692842619e-08) (12, 3.75807605603847e-08) (15, 7.54601925261432e-08) (18, 1.32463724821618e-07) (21, 2.12905975498591e-07) (24, 3.20571371573657e-07) (27, 4.59065146736666e-07) (30, 6.33248126619824e-07) (33, 8.46355638687303e-07) (36, 1.10293157122632e-06) (39, 1.40628134320465e-06) (42, 1.76111852522038e-06) (45, 2.17112088817838e-06) (48, 2.64001232856044e-06)};
\legend {(a), (b), (c), (d), (e), (f), (g), (h)}
\end{axis}
\end{tikzpicture}
\end{subfigure}
\begin{subfigure}[b]{0.47\textwidth}
\begin{tikzpicture}
\begin{axis}[
    xlabel={$\lambda$}, ylabel={$\Delta_{rel}$}, ymode = log, legend pos= south east ]
\addplot coordinates { (2, 0.00338851406217308) (4, 0.00717399746424729) (6, 0.00782557626621457) (8, 0.013891754269775) (9, 0.0106820281722022)};
\addplot coordinates { (2, 0.000830222022899685) (4, 0.00266611235997507) (6, 0.00229300076494109) (9, 0.00351399553783415) (12, 0.00586110160056273) (15, 0.00904173741630135) (16, 0.0102556301584657)};
\addplot coordinates { (2, 0.000562691069476195) (4, 0.00122585315636946) (6, 0.00175698507367345) (9, 0.00371530978089724) (12, 0.00689821955727446) (15, 0.0110029103691482) (18, 0.0158914680543967) (21, 0.0215325674581223) (24, 0.0280951885245139) (27, 0.0361451724204803) (30, 0.0469959493934096) (33, 0.0633685503999436) (35, 0.0798890059674621)};
\addplot coordinates { (2, 0.000533247519255101) (4, 0.00098843553097078) (6, 0.00161445825372413) (9, 0.00368396925410514) (12, 0.00688998419442918) (15, 0.0110015246206132) (18, 0.0158914145478577) (21, 0.0215325737009227) (24, 0.0280953442025205) (27, 0.03614544219651) (30, 0.0469960830435601) (33, 0.0633685564351415) (36, 0.090798081039622) (39, 0.141237763411618) (42, 0.245009553303777) (45, 0.510778439515317) (46, 0.708214137184029) (48, 1.96028495276276)};
\addplot coordinates { (2, 4.99376563684731e-05) (4, 0.000121970317326893) (6, 0.000140032806289974) (9, 0.000249082038604912) (12, 0.00043698019938019) (15, 0.000682753721352085) (18, 0.000982384529157063) (21, 0.00115053429061055) (24, 0.00145725528008062) (27, 0.00180543645497172) (30, 0.00218496687244472) (33, 0.00318536845268981) (36, 0.00302878173362933) (39, 0.00443978512347803) (42, 0.00396696639976743) (45, 0.00589739617984154) (48, 0.00670668612897734)};
\addplot coordinates { (2, 5.01542245845209e-07) (4, 1.15356679288429e-06) (6, 1.37529287864765e-06) (9, 2.36824814033448e-06) (12, 4.07023976799231e-06) (15, 6.24496390663589e-06) (18, 8.83519837596856e-06) (21, 1.18317121325434e-05) (24, 1.52349988538938e-05) (27, 1.90475102038972e-05) (30, 2.32716195220911e-05) (33, 2.79090266537701e-05) (36, 3.2960297160293e-05) (39, 3.84279436398156e-05) (42, 4.43185020257375e-05) (45, 5.06480566329206e-05) (48, 5.74304997162916e-05)};
\addplot coordinates { (2, 4.86951710751613e-09) (4, 1.17103789806649e-08) (6, 1.35734518459633e-08) (9, 2.38842914053757e-08) (12, 4.1235073327516e-08) (15, 6.34262441370206e-08) (18, 8.99358479268478e-08) (21, 1.20692051413573e-07) (24, 1.55716248784086e-07) (27, 1.95048809718894e-07) (30, 2.38729378552588e-07) (33, 2.86791231752452e-07) (36, 3.3887910731306e-07) (39, 3.97217399793917e-07) (42, 4.68522670075778e-07) (45, 5.69375535714813e-07) (48, 7.14360970602078e-07)};
\addplot coordinates { (2, 4.90322168273426e-11) (4, 1.16915597808861e-10) (6, 1.3593141841574e-10) (9, 2.38745015710512e-10) (12, 4.12100868216076e-10) (15, 6.33796991320739e-10) (18, 8.98608062423592e-10) (21, 1.20581815820193e-09) (24, 1.55564666012857e-09) (27, 1.94850408403659e-09) (30, 2.38479504381334e-09) (33, 2.86489781825019e-09) (36, 3.00760880466253e-09) (39, 5.01709633421235e-09) (42, 4.56841463215627e-09) (45, 5.24659308421796e-09) (48, 5.86484887681805e-09)};
\end{axis}
\end{tikzpicture}
\end{subfigure}
\caption{The Troesch's problem. Relative errors of approximations for $u_{2}(0)$ (to the left) and $u_{2}(1)$ (to the right) calculated using $I$-$\mathcal{SP}_{1}\mathcal{FP}_{2}$ transformation strategy and the mesh refinement procedure with $M = 0.1$ and different values of $h_{max}$ and $h_{min}$ (see Tab.~\ref{config_table} to get the description of the legend abbreviations). } 
\label{Example_1_graph_delta_transformation_old}
\end{figure}
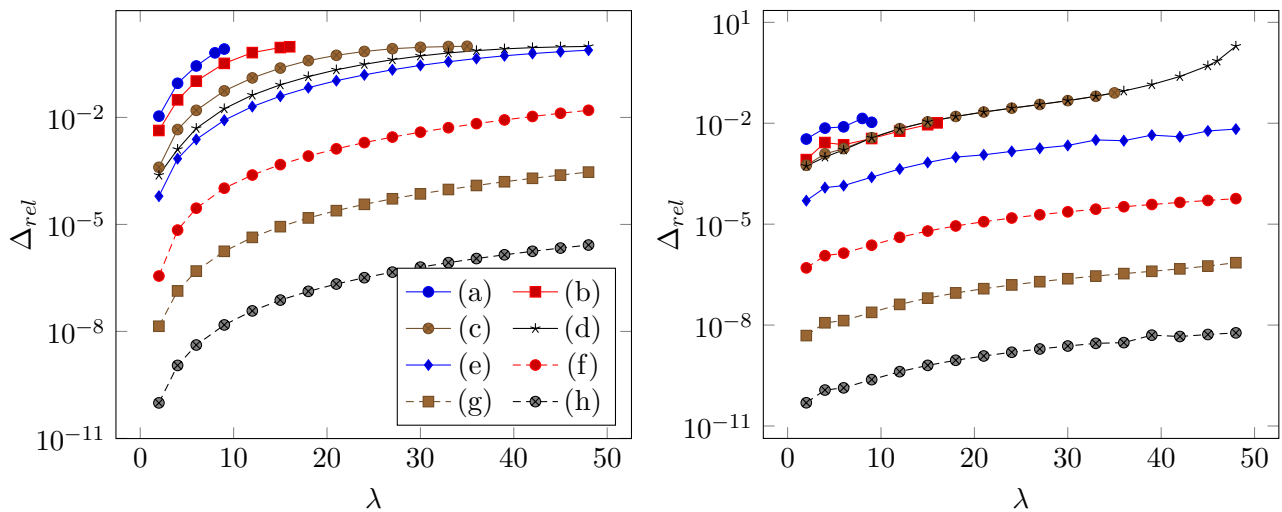

\begin{table}[]
\begin{tabular}{|l|l|l|l|l|l|l|l|l|}
\hline
\backslashbox{{\small Parameter}}{{\small Configuration}} & (a) & (b) & (c) & (d) & (e) & (f) & (g) & (h) \\
\hline
$h_{min}$       & 0.1 & 0.05 & 0.025 & 0.01 & 0.01 & $10^{-3}$ & $10^{-4}$ & $10^{-5}$\\
$h_{max}$       & 0.1 & 0.1  & 0.1   & 0.1  & 0.01 & $10^{-3}$ & $10^{-4}$ & $10^{-5}$\\
Max. mesh size  & 15  &  28  & 57    & 119  & 147  & 1489      & 12.6K     & 132K     \\
 \hline
\end{tabular}
\caption{ Characteristics of different mesh refinement configurations mentioned in the legend of Fig.~\ref{Example_1_graph_delta_transformation_old}. } \label{config_table}
\end{table}

The relative errors of approximations for values $u_{2}(0)$ and $u_{2}(1)$ obtained with $I$-$\mathcal{SP}_{1}\mathcal{FP}_{2}$ transformation strategy are shown in Fig.~\ref{Example_1_graph_delta_transformation_old}. Notice that each error curve in the figure correspond to certain values of parameters $h_{min},$ $h_{max}$ (used in the mesh refinement procedure) and those can be derived from Tab.~\ref{config_table}. The first thing that can be pointed out, when looking onto the corresponding charts, is that when using strategy $I$-$\mathcal{SP}_{1}\mathcal{FP}_{2}$ we are still able to get approximation errors of order $\mathcal{O}(h^{2}),$ just as it is for the regular trapezoidal scheme. Indeed, 10 times increase in the number of mesh points (see Tab.~\ref{config_table}) results in 100 times decline of the approximation error magnitudes. Another important observation to make is that the error curves in Fig.~\ref{Example_1_graph_delta_transformation_old} are sloped down to a considerably greater extent than similar curves corresponding to the regular trapezoidal scheme (see Fig.~\ref{Example_1_graph_delta_uniform}, \ref{Example_1_graph_delta_nonuniform}). This clearly indicates that      strategy $I$-$\mathcal{SP}_{1}\mathcal{FP}_{2}$ is considerably less sensitive to the increase of parameter $\lambda$ (i.e., possesses a higher level of "stiffness resistance").

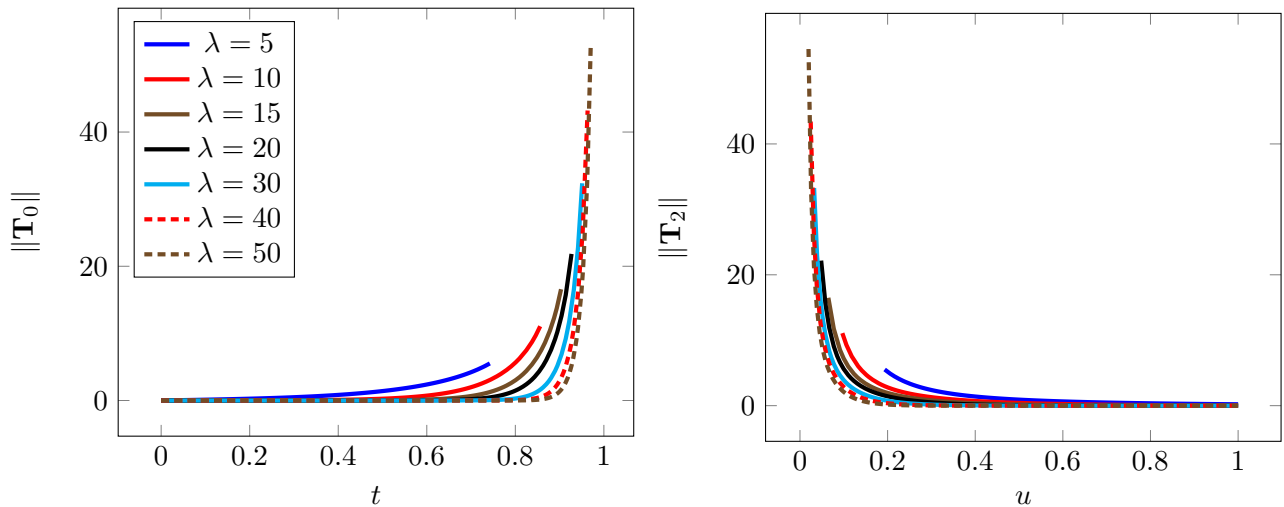
\begin{figure}[h!]
\centering
\begin{subfigure}[b]{0.47\textwidth}
\begin{tikzpicture}
\begin{axis}[no markers,
    every axis plot/.append style={ultra thick}, xlabel={$t$}, ylabel={$\|\mathbf{T}_{0}\|$}, legend pos= north west]
\addplot table [x index = {4}, y index = {5}, col sep=comma] {Troesch_rhs_lambda_5.csv};
\addplot table [x index = {4}, y index = {5}, col sep=comma] {Troesch_rhs_lambda_10.csv};
\addplot table [x index = {4}, y index = {5}, col sep=comma] {Troesch_rhs_lambda_15.csv};
\addplot table [x index = {4}, y index = {5}, col sep=comma] {Troesch_rhs_lambda_20.csv};
\addplot[cyan] table [x index = {4}, y index = {5}, col sep=comma] {Troesch_rhs_lambda_30.csv};
\addplot table [x index = {4}, y index = {5}, col sep=comma] {Troesch_rhs_lambda_40.csv};
\addplot table [x index = {4}, y index = {5}, col sep=comma] {Troesch_rhs_lambda_50.csv};
\legend { $\lambda = 5$, $\lambda = 10$, $\lambda = 15$, $\lambda = 20$, $\lambda = 30$, $\lambda = 40$, $\lambda = 50$}
\end{axis}
\end{tikzpicture}
\end{subfigure}
\begin{subfigure}[b]{0.47\textwidth}
\begin{tikzpicture}
\begin{axis}[no markers,
    every axis plot/.append style={ultra thick}, xlabel={$u$}, ylabel={$\|\mathbf{T}_{2}\|$}, legend pos= south east]
\addplot table [x index = {2}, y index = {3}, col sep=comma] {Troesch_rhs_lambda_5.csv};
\addplot table [x index = {2}, y index = {3}, col sep=comma] {Troesch_rhs_lambda_10.csv};
\addplot table [x index = {2}, y index = {3}, col sep=comma] {Troesch_rhs_lambda_15.csv};
\addplot table [x index = {2}, y index = {3}, col sep=comma] {Troesch_rhs_lambda_20.csv};
\addplot[cyan] table [x index = {2}, y index = {3}, col sep=comma] {Troesch_rhs_lambda_30.csv};
\addplot table [x index = {2}, y index = {3}, col sep=comma] {Troesch_rhs_lambda_40.csv};
\addplot table [x index = {2}, y index = {3}, col sep=comma] {Troesch_rhs_lambda_50.csv};
\end{axis}
\end{tikzpicture}
\end{subfigure}
\caption{The Troesch's problem. $I$-$\mathcal{SP}_{1}\mathcal{FP}_{2}$ transformation strategy. Norms of the right-hand side vector-functions $\mathbf{T}_{0} = \mathbf{T}_{0}(u_{1}(t), u_{2}(t), \lambda)$ (to the left, see \eqref{Troesch_problem}, \eqref{Troesch_problem_left_part_v1}) and $\mathbf{T}_{2} = \mathbf{T}_{2}(w_{1}(u), w_{2}(u), \lambda)$ (to the right, see \eqref{Troesch_problem_right_part_v3}), as functions of $t$ and $u$ respectively, evaluated for different values of $\lambda$.} 
\label{Example_1_graph_rhs_old}
\end{figure}

It is easy to see that, in the vicinity of point $t=1,$ the norm of vector-function $\mathbf{T}(u_{1}(t), u_{2}(t), \lambda)$ \eqref{Troesch_problem} is an exponentially increasing function of $\lambda$ (provided that functions $u_{1}(t)$ and $u_{2}(t)$ denote the solution to the Troesch's problem). This is not the case for the right-hand side vector-functions of the transformed problem \eqref{Troesch_problem_left_part_v1}, \eqref{Troesch_problem_right_part_v3}, \eqref{Matching_and_boundary_conditions_Troesch_v3}. As it can be seen in Fig.~\ref{Example_1_graph_rhs_old}, the corresponding norms take their maximum values exactly at the matching point $\varepsilon$ \eqref{first_matching_point_condition} and their magnitudes are comparable to the corresponding values of parameter $\lambda.$ According to the definition of stiffness accepted throughout the present paper, the latter confirms that the transformed problem is considerably less stiff than the initial one. However, the "bump" near the matching point steadily "grows" along with $\lambda$ (see Fig.~\ref{Example_1_graph_rhs_old}). Although the refinement procedure can mitigate this issue by putting more mesh points to that region, in theory, for extremely high values of parameter $\lambda,$ this can cause problems to the underlying finite difference scheme. The issue would have been resolved if we had used transformation $\mathcal{SP}_{2}$ instead of $\mathcal{SP}_{1}\mathcal{FP}_{2}.$ Unfortunately, as it was pointed out in section \ref{transformations_section}, transformation $\mathcal{SP}_{2},$ if applied all the way to the point $t=1,$ would result in a BVP with a non-stationary boundary point which we agreed to leave for the future studies. Nevertheless, there is a way for us to take advantage of $\mathcal{SP}_{2}$ and this is by using yet another transformation strategy which can be abbreviated as $I$-$\mathcal{SP}_{2}$-$\mathcal{SP}_{1}\mathcal{FP}_{2}.$ 

\subsubsection{Strategy $I$-$\mathcal{SP}_{2}$-$\mathcal{SP}_{1}\mathcal{FP}_{2}$}
As one could deduce, strategy $I$-$\mathcal{SP}_{2}$-$\mathcal{SP}_{1}\mathcal{FP}_{2}$ implies introducing one more sub-interval, $[\varepsilon, \varepsilon_{1}] \subset (0, 1),$ on which transformation $\mathcal{SP}_{2}$ is to be applied. The overall result of the new strategy amounts to the following BVP:
\begin{eqnarray}
\left[\begin{array}{ll}
u_{1}^{\prime}(t), u_{2}^{\prime}(t)
\end{array}\right]^{T} & = & \mathbf{T}_{0}(u_{1}(t), u_{2}(t), \lambda), \; t\in [0, \varepsilon], \; \text{(identity transformation part)}; \nonumber \\ 
\left[\begin{array}{ll}
v_{1}^{\prime}(\mu), v_{2}^{\prime}(\mu)
\end{array}\right]^{T} & = & \mathbf{T}_{1}(v_{1}(\mu), v_{2}(\mu), \lambda), \; \mu \in [u_{2}(\varepsilon), u_{2}(\varepsilon_{1})], \; \text{($\mathcal{SP}_{2}$ transformation part)}; \label{Troesch_transformed_new_strategy} \\
\left[\begin{array}{ll}
w_{1}^{\prime}(\nu), w_{2}^{\prime}(\nu)
\end{array}\right]^{T} & = & \mathbf{T}_{2}(w_{1}(\nu), w_{2}(\nu), \lambda), \; \nu \in [u_{1}(\varepsilon_{1}), 1], \; \text{($\mathcal{SP}_{1}\mathcal{FP}_{2}$ transformation part)}; \nonumber
\end{eqnarray}
$$v_{1}(u_{2}(\varepsilon)) = u_{1}(\varepsilon); \; v_{2}(u_{2}(\varepsilon))= \varepsilon, \;\;\; v_{1}(u_{2}(\varepsilon_{1})) = u_{1}(\varepsilon_{1}), \; v_{2}(u_{2}(\varepsilon_{1}))= \varepsilon_{1};$$
$$w_{2}(u_{1}(\varepsilon_{1})) = \frac{1}{u_{2}(\varepsilon)}, \; w_{1}(u_{1}(\varepsilon_{1}))= \varepsilon_{1};$$
$$u_{1}(0) = 0,\; w_{1}(1) = 1,$$
where vector-functions $\mathbf{T}_{i}(\cdot), \; i =0,1,2$ have been introduced in \eqref{Troesch_problem}, \eqref{Troesch_problem_right_part_v1} and \eqref{Troesch_problem_right_part_v3} respectively.

\begin{figure}[h!]
\centering
\begin{subfigure}[b]{0.3\textwidth}
\begin{tikzpicture}
\begin{axis}[width=6.25cm,height=5cm, no markers,
    every axis plot/.append style={ultra thick}, 
    y label style={at={(axis description cs:0.3,1.0)}, rotate=-90, anchor=south},
    xlabel={$t$}, ylabel={$\|\mathbf{T}_{0}\|$}, legend pos= north west]
\addplot table [x index = {6}, y index = {7}, col sep=comma] {Troesch_new_rhs_lambda_10.csv};
\addplot table [x index = {6}, y index = {7}, col sep=comma] {Troesch_new_rhs_lambda_20.csv};
\addplot table [x index = {6}, y index = {7}, col sep=comma] {Troesch_new_rhs_lambda_30.csv};
\addplot table [x index = {6}, y index = {7}, col sep=comma] {Troesch_new_rhs_lambda_40.csv};
\addplot[cyan] table [x index = {6}, y index = {7}, col sep=comma] {Troesch_new_rhs_lambda_50.csv};
\legend { $\lambda = 10$, $\lambda = 20$, $\lambda = 30$, $\lambda = 40$, $\lambda = 50$}
\end{axis}
\end{tikzpicture}
\end{subfigure}
\begin{subfigure}[b]{0.3\textwidth}
\begin{tikzpicture}
\begin{axis}[width=6.25cm,height=5cm, no markers,
    every axis plot/.append style={ultra thick}, 
    y label style={at={(axis description cs:0.3,1.0)}, rotate=-90, anchor=south},  xmode = log,
    xlabel={$\mu$}, ylabel={$\|\mathbf{T}_{1}\|$}, legend pos= south east]
\addplot table [x index = {4}, y index = {5}, col sep=comma] {Troesch_new_rhs_lambda_10.csv};
\addplot table [x index = {4}, y index = {5}, col sep=comma] {Troesch_new_rhs_lambda_20.csv};
\addplot table [x index = {4}, y index = {5}, col sep=comma] {Troesch_new_rhs_lambda_30.csv};
\addplot table [x index = {4}, y index = {5}, col sep=comma] {Troesch_new_rhs_lambda_40.csv};
\addplot[cyan] table [x index = {4}, y index = {5}, col sep=comma] {Troesch_new_rhs_lambda_50.csv};
\end{axis}
\end{tikzpicture}
\end{subfigure}
\begin{subfigure}[b]{0.3\textwidth}
\begin{tikzpicture}
\begin{axis}[width=6.25cm,height=5cm, no markers,
    every axis plot/.append style={ultra thick}, 
    y label style={at={(axis description cs:0.3,1.0)}, rotate=-90, anchor=south},
    xlabel={$\nu$}, ylabel={$\|\mathbf{T}_{2}\|$}, legend pos= south east]
\addplot table [x index = {2}, y index = {3}, col sep=comma] {Troesch_new_rhs_lambda_10.csv};
\addplot table [x index = {2}, y index = {3}, col sep=comma] {Troesch_new_rhs_lambda_20.csv};
\addplot table [x index = {2}, y index = {3}, col sep=comma] {Troesch_new_rhs_lambda_30.csv};
\addplot table [x index = {2}, y index = {3}, col sep=comma] {Troesch_new_rhs_lambda_40.csv};
\addplot[cyan] table [x index = {2}, y index = {3}, col sep=comma] {Troesch_new_rhs_lambda_50.csv};
\end{axis}
\end{tikzpicture}
\end{subfigure}
\caption{The Troesch's problem. $I$-$\mathcal{SP}_{2}$-$\mathcal{SP}_{1}\mathcal{FP}_{2}$ transformation strategy. Norms of the right-hand side vector-functions $\mathbf{T}_{0} = \mathbf{T}_{0}(u_{1}(t), u_{2}(t), \lambda)$ (to the left, see \eqref{Troesch_problem}), $\mathbf{T}_{1} = \mathbf{T}_{1}(v_{1}(\mu), v_{2}(\mu), \lambda)$ (in the middle, see \eqref{Troesch_problem_right_part_v1}) and $\mathbf{T}_{2} = \mathbf{T}_{2}(w_{1}(\nu), w_{2}(\nu), \lambda)$ (to the right, see \eqref{Troesch_problem_right_part_v3}), as functions of $t,$ $\mu$ and $\nu$ respectively, evaluated for different values of $\lambda$.} 
\label{Example_1_graph_rhs_new}
\end{figure}
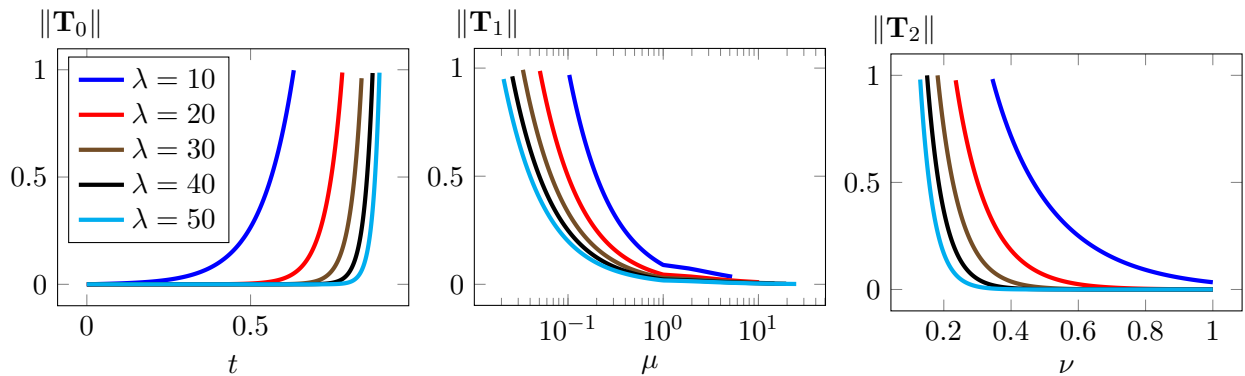

\begin{figure}[h!]
\centering
\begin{subfigure}[b]{0.47\textwidth}
\begin{tikzpicture}
\begin{axis}[ 
    xlabel={$\lambda$}, ylabel={$\Delta_{rel}$}, ymode = log, legend pos= south east 
    ]
\addplot coordinates { (12, 4.66183808576599e-05) (13, 4.92277331843321e-05) (15, 5.19736651339117e-05) (18, 4.10756801268261e-05) (20, 2.50173686601479e-05) (21, 1.2757014411789e-05) (22, 5.55661789907544e-06) (23, 2.63508021432888e-05) (24, 5.23202616270393e-05) (27, 0.000150436931891877) (30, 0.000302457933859424) (33, 0.000500976852219944) (36, 0.000766582194330814) (39, 0.00108791499140681) (42, 0.00149053576489454) (45, 0.00196099025714519) (48, 0.00250585505103846)};
\addplot coordinates { (12, 6.56238329295924e-05) (13, 8.42741089712588e-05) (15, 0.000131910002796506) (18, 0.000232646509844993) (20, 0.000322865062255045) (21, 0.000375257206884869) (22, 0.000432960217230155) (23, 0.000496059383676925) (24, 0.000564922017529833) (27, 0.000810133680084343) (30, 0.00111392146968859) (33, 0.0014901961475156) (36, 0.00194367791116476) (39, 0.00247682432154792) (42, 0.0030937911286009) (45, 0.00381616172211744) (48, 0.0046435074216597)};
\legend {$I$-$\mathcal{SP}_{2}$-$\mathcal{SP}_{1}\mathcal{FP}_{2}$, $I$-$\mathcal{SP}_{1}\mathcal{FP}_{2}$}
\end{axis}
\end{tikzpicture}
\end{subfigure}
\caption{The Troesch's problem. Relative error of the approximations for $u_{2}(0)$ calculated for different values of $\lambda$ using $I$-$\mathcal{SP}_{1}\mathcal{FP}_{2}$ and $I$-$\mathcal{SP}_{2}$-$\mathcal{SP}_{1}\mathcal{FP}_{2}$ transformation strategies on the same meshes (without mesh refinement). } 
\label{Example_1_comparison_of_strategies}
\end{figure}
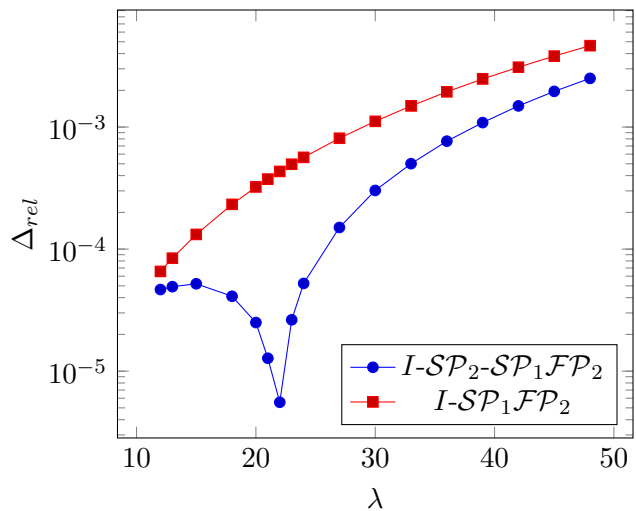

By choosing $\varepsilon$ and $\varepsilon_{1}$ as follows
$$\varepsilon = \min\limits_{t\in [0,1]}\left\{t : \lambda \sinh\left(\lambda u_{1}(t)\right) > 1\right\},$$ 
$$\varepsilon_{1} = \min\limits_{t\in [0,1]}\left\{t : \lambda \sinh\left(\lambda u_{1}(t)\right) < \left(u_{2}(t)\right)^{3}\right\},$$
we, quite predictably, get the results shown in Fig.~\ref{Example_1_graph_rhs_new}. Obviously, with a proper choice of $\varepsilon$ and $\varepsilon_{1},$ strategy $I$-$\mathcal{SP}_{2}$-$\mathcal{SP}_{1}\mathcal{FP}_{2}$ is able to keep the norms of the right-hand sides of equations \eqref{Troesch_transformed_new_strategy} "below" 1, for any $\lambda > 0.$ Comparing the graphs in  Fig.~\ref{Example_1_graph_rhs_old} and \ref{Example_1_graph_rhs_new}, one can argue that the "bump" around the matching point of $I$-$\mathcal{SP}_{1}\mathcal{FP}_{2}$ was eliminated by means of stretching the integration interval. The latter can be clearly seen on the middle chart in Fig.~\ref{Example_1_graph_rhs_new}, having its horizontal axis scaled logarithmically. Thus, although problem \eqref{Troesch_transformed_new_strategy} is "less stiff" than the one strategy $I$-$\mathcal{SP}_{1}\mathcal{FP}_{2}$ results in (according to the definition of stiffness adopted in this paper), it is still not clear whether we really gain any benefits (from the computational standpoint) when using $I$-$\mathcal{SP}_{2}$-$\mathcal{SP}_{1}\mathcal{FP}_{2}.$ This was not investigated in depth and a thorough comparison of the two strategies is left for the future publications. Based on a number of experiments, we can conclude that the naive mesh refinement procedure used in this section is definitely sub-optimal when applied to problem \eqref{Troesch_transformed_new_strategy} (it puts unnecessarily many points on the interval $[u_{2}(\varepsilon), u_{2}(\varepsilon_{1})]$ associated with the second equation). The question about more optimal mesh refinement algorithm, that can be used in this case, is still open as well as the question about optimal choice of the matching points $\varepsilon$ and $\varepsilon_{1}.$ At the same time, numerical experiments also suggest that strategy $I$-$\mathcal{SP}_{2}$-$\mathcal{SP}_{1}\mathcal{FP}_{2}$ provides lower approximation errors than strategy $I$-$\mathcal{SP}_{1}\mathcal{FP}_{2},$ when both are applied on the same meshes (i.e. when the mesh is fixed and no refinement is used). The latter is especially pronounced near the point $t=0,$ which is illustrated in Fig. \ref{Example_1_comparison_of_strategies}.

In conclusion to this section, as it was promised earlier, we provide some results obtained with a multi-precision arithmetic, namely, the one based on the floating-point type \verb"number<cpp_dec_float<110>, et_off>" from Boost C++ Library (see \url{https://www.boost.org}). The template argument, "110", determines the decimal digits precision provided by the type (the machine epsilon in this case turns out to be of order $10^{-130}$). The corresponding approximations of $u_{2}(0)$ and $u_{2}(1)$ are presented in Tab.~\ref{data_at_0_table} and \ref{data_at_1_table} respectively. The calculations are conducted on a series of meshes with different discretization levels: each next mesh contains approximately 10 times more points than the previous one (see the leftmost columns of the two tables). Each column of tables \ref{data_at_0_table} and \ref{data_at_1_table} (except for the leftmost ones) contains numerical data calculated for a certain value of parameter $\lambda,$ which is specified in the column's header together with a multiplication factor that needs to be applied to the data values in order to get the actual approximations for $u_{2}(0)$ and $u_{2}(1).$ Colouring is used to visualize "trusted" digits of the approximations. The number of leading green digits in each numerical value is equal to the negative order of the corresponding relative error evaluated with respect to the reference data (which, in this case, is the data from \cite{Gen_sol_of_TP}, see the bottommost row in Tab.~\ref{data_at_0_table}). For example, $k$ leading green digits in a value means that the corresponding relative error is of order $10^{-k}.$ Cyan color serves the same purpose with the only difference that it is used when reference data is not available, in which case the relative error was estimated indirectly based on the convergence considerations (evaluating differences between approximations obtained on successive meshes) and an assumption that the approximation error of $I$-$\mathcal{SP}_{2}$-$\mathcal{SP}_{1}\mathcal{FP}_{2}$ behaves as $\mathcal{O}(h^{2})$ (i.e., a mesh with 10 times more points, should result in 2 more "trusted" digits in the corresponding approximation).

\begin{table}[]
\begin{tabular}{|l|l|l|l|l|l|l|l|l|}
\hline
Mesh  & $\lambda = 50;$ & $\lambda = 100;$ & $\lambda = 200;$ & $\lambda = 300;$ & $\lambda = 400;$ & $\lambda = 500;$  \\
size  & $\times 10^{-21}$ & $\times 10^{-43}$ & $\times 10^{-86}$ & $\times 10^{-130}$ &  $\times 10^{-173}$ & $\times 10^{-217}$  \\
\hline
1.4K       & \textcolor{green}{1.4}93019771 & \textcolor{green}{2}.283685233 & \textcolor{green}{1}.280634173 & --- & --- &--- \\
10K        & \textcolor{green}{1.542}514970 & \textcolor{green}{2.96}7678751 & \textcolor{green}{1.0}80997285 & \textcolor{cyan}{3.7}945093 & \textcolor{cyan}{1.2}595406 & \textcolor{green}{3}.880275407\\
100K       & \textcolor{green}{1.54299}5447 & \textcolor{green}{2.9759}78741 & \textcolor{green}{1.106}851674 & \textcolor{cyan}{4.115}1391 & \textcolor{cyan}{1.529}0852 & \textcolor{green}{5.67}7350512\\
1M         & \textcolor{green}{1.5429998}34 & \textcolor{green}{2.976059}973 & \textcolor{green}{1.10711}4604 & \textcolor{cyan}{4.11852}62 & \textcolor{cyan}{1.53210}53 & \textcolor{green}{5.6994}38135\\
10M        & \textcolor{green}{1.542999878} & \textcolor{green}{2.97606077}3 & \textcolor{green}{1.1071171}95 & \textcolor{cyan}{4.1185598} & \textcolor{cyan}{1.5321354} & \textcolor{green}{5.699658}911\\
 \hline
 \hline
\cite{Gen_sol_of_TP} &  1.542999878 & 2.976060781 & 1.107117221 & --- & --- & 5.699661125 \\
 \hline
\end{tabular}
\caption{ Troesch's problem. Approximations of $u_{2}(0)$ for some values of $\lambda$ calculated with $I$-$\mathcal{SP}_{2}$-$\mathcal{SP}_{1}\mathcal{FP}_{2}$ transformation strategy on different meshes using multiple-precision arithmetic. } \label{data_at_0_table}
\end{table}

\begin{table}[]
\begin{tabular}{|l|l|l|l|l|l|l|l|l|}
\hline
Mesh  & $\lambda = 50;$ & $\lambda = 100;$ & $\lambda = 200;$ & $\lambda = 300;$ & $\lambda = 400;$ & $\lambda = 500;$  \\
size  & $\times 10^{10}$ & $\times 10^{21}$ & $\times 10^{43}$ & $\times 10^{65}$ &  $\times 10^{86}$ & $\times 10^{108}$  \\
\hline
1.4K       & \textcolor{cyan}{7.2009}32187365 & \textcolor{cyan}{5.186}00902647 & \textcolor{cyan}{2.69}085262858 & \textcolor{cyan}{1.396}90969485 & \textcolor{cyan}{7.25}54758236 & \textcolor{cyan}{3.77}03273898 \\
10K        & \textcolor{cyan}{7.200496}850882 & \textcolor{cyan}{5.18472}589782 & \textcolor{cyan}{2.68815}990996 & \textcolor{cyan}{1.39375}966426 & \textcolor{cyan}{7.2264}361921 & \textcolor{cyan}{3.7468}295566 \\
100K       & \textcolor{cyan}{7.20049000}3218 & \textcolor{cyan}{5.1847057}3232 & \textcolor{cyan}{2.6881175}6950 & \textcolor{cyan}{1.3937100}8152 & \textcolor{cyan}{7.225978}3926 & \textcolor{cyan}{3.746458}3642 \\
1M         & \textcolor{cyan}{7.2004899344}33 & \textcolor{cyan}{5.184705530}62 & \textcolor{cyan}{2.688117146}09 & \textcolor{cyan}{1.393709585}67 & \textcolor{cyan}{7.22597381}44 & \textcolor{cyan}{3.74645465}20 \\
10M        & \textcolor{cyan}{7.200489933746} & \textcolor{cyan}{5.18470552861} & \textcolor{cyan}{2.68811714186} & \textcolor{cyan}{1.39370958072} & \textcolor{cyan}{7.2259737686} & \textcolor{cyan}{3.7464546149} \\
 \hline
\end{tabular}
\caption{ Troesch's problem. Approximations of $u_{2}(1)$ for some values of $\lambda$ calculated with $I$-$\mathcal{SP}_{2}$-$\mathcal{SP}_{1}\mathcal{FP}_{2}$ transformation strategy on different meshes using multiple-precision arithmetic. } \label{data_at_1_table}
\end{table}

\section{Conclusions}\label{Conclusions_section}
The results of the numerical examples in the previous section clearly show that the transformation-based approach can serve as a useful and quite powerful tool for solving stiff two-point boundary value problems of a certain type. Potentially, the suggested methodology can be used in conjunction with absolutely any existing numerical method for solving two-point BVPs, helping to overcome computational difficulties associated with the boundary and shock layers. 

Together with some answers given in the present paper there are also a few questions posed in it. Those can be seen as landmarks suggesting the directions in which the theory can be developed further. For instance, it seems quite attractive (a) to investigate the prospects of using transformations that result in BVPs with non-stationary boundary points (see section \ref{transformations_section}), (b) to study the strengths and weaknesses of different transformation strategies, (c) to design efficient/optimal mesh refinement algorithms for each particular strategy. Looking forward, it also makes a good sense to investigate a possibility of extending the suggested methodology to the problems with non-local boundary conditions as well as to the systems of partial differential equations.

\bibliography{../bibliography/references_stat}{}
\bibliographystyle{plain}
\end{document}